\documentclass{amsart}
\usepackage{amsfonts,amscd}

\newtheorem{claim}{}[section]
\newtheorem{theorem}[claim]{Theorem}
\newtheorem{lemma}[claim]{Lemma}
\newtheorem{proposition}[claim]{Proposition}
\newtheorem{corollary}[claim]{Corollary}

\def\proclaim #1. #2\par{\medbreak
\noindent{\bf#1.\enspace}{\sl#2}\par\medbreak}
\makeatother

\begin{document}
\title[One-sided $M$-Ideals and Multipliers in Operator Spaces]{One-sided
$M$-Ideals and
Multipliers \\ in Operator Spaces, I}
\thanks{This is a small modification of the April 2001 revision.  It is not the
final published version}
\author{David P. Blecher}
\address{Department of Mathematics, University of Houston, Houston,
TX 77204-3476}
\email[David P. Blecher]{dblecher@math.uh.edu}
\author{Edward G. Effros}
\address{Department of Mathematics, UCLA, Los Angeles, CA 90095-1555}
\email[Edward G. Effros]{ege@math.ucla.edu}
\author{Vrej Zarikian}
\address{Department of Mathematics,
University of Texas at Austin, Austin, TX 78712-1082}
\email[Vrej Zarikian]{zarikian@math.utexas.edu}
\thanks{Blecher and Effros were partially supported by the National Science
Foundation}
\begin{abstract}
The theory of $M$-ideals and multiplier mappings of Banach spaces
naturally generalizes to
left (or right) $M$-ideals and multiplier mappings of operator
spaces. These subspaces and
mappings are intrinsically characterized in terms of the matrix
norms. In turn this is used
to prove that the algebra of left adjointable mappings of a dual
operator space $X$ is a
von Neumann algebra. If in addition $X$ is an operator
$A$--$B$-bimodule for $C^{*}$-algebras $A$ and $B$,
then the module operations
on $X$ are automatically weak$^{*}$ continuous. One sided $L$-projections are
introduced, and analogues of various results from the classical
theory are proved. An
assortment of examples is considered.
\end{abstract}
\maketitle

\let\text=\mbox

\section{Introduction}

\label{Int}

It has long been recognized that the algebraic structure of a $C^{*}$%
-algebra $A$ is closely linked to its geometry as a Banach space (see
\cite{Kadison}). This principle was illustrated in \cite{AlfsenEffros},
and \cite{Alfsen}, p. 237,
where it was shown that the closed two-sided ideals of a
$C^{*}$-algebra
coincide with the $M$\emph{-ideals} of the underlying Banach space
(see also \cite{SmithWard}).
Similarly, the center of a $C^{*}$-algebra is determined by the
\emph{centralizer mappings} of the Banach space \cite{AlfsenEffros},
\cite{Behrends}. It was
subsequently shown that these notions can be applied to a broad range of
Banach space problems unrelated to operator algebra theory (see
\cite{HWW} for
references to the extensive literature on this subject).

In this paper we show that one can similarly characterize the closed
one-sided ideals and one-sided multipliers in a $C^{*}$-algebra in terms of
its \emph{matrix norms}, i.e. its underlying \emph{operator space
structure}.
We show that
the closed one-sided ideals in a $C^{*}$-algebra are just the
\emph{complete
one-sided} $M$-\emph{ideals }(defined below) of the operator space. We
also prove that the \emph{one-sided multipliers} and the \emph{one-sided
adjointable multipliers} of an operator space (first studied independently in
\cite{BlecherShilov} and \cite{WWerner}, see also \cite{BlecherPaulsen})
    have surprisingly simple matrix norm characterizations.
Once again these abstract considerations have
important applications elsewhere, including a striking automatic
continuity result for
dual modules (see Corollary \ref{automcont}). They have also led to a new
characterization of the dual operator algebras {\cite{Blecherdual}.

Turning to the details, if $X$ is an operator space, a linear mapping
$P:X\rightarrow X$ with $P^{2}=P$ is said to be a \emph{left
}$M$\emph{-projection}
if for each $x\in X,$
\begin{equation*}
\left\| x\right\| =\left\| \left[
\begin{array}{c}
P(x) \\
x-P(x)
\end{array}
\right] \right\| .
\end{equation*}
We say that
$P$ is a \emph{complete} left $M$-projection if for each $n\in \Bbb{N},$ $
P_{n}:M_{n}(X)\rightarrow M_{n}(X)$ is a left $M$-projection.
Here $P_n$ is the canonical ``entry-wise'' action of $P$ on matrices.
A subspace $J$
of $X$ is a (complete) \emph{right }$M$-\emph{summand} if $J=P(X)$ with $P$
a (complete) left $M$-projection. Finally, a closed subspace $J$ of $X$ is a
(complete) \emph{right} $M$-\emph{ideal} if $J^{\perp \perp }$ is a
(complete) right $M$-summand. If $A$ is a unital $C^{*}$-algebra,
then the complete left $M$-projections
are given by $P(x)=ex$ where $e$ is an orthogonal  projection in $A.$
Hence the
complete right $M$-summands of $A$ are
the algebraic right ideals of the form $eA.$
As a consequence the complete right $M$-ideals in a
$C^{*}$-algebra are exactly the closed right ideals.
One may similarly define the notion of a
{\em right $M$-projection} by using row matrices. We have left the routine
details of such reversed notions (left $M$-summands, etc.) to the reader.

As in the theory of $M$-ideals in a Banach space, it is technically useful
to introduce the dual notions of \emph{one-sided} \emph{$L$-projections},
\emph{$L$-summands} and
$L$\emph{-ideals }in an operator space. We also prove
that complete one-sided
$L$-ideals are necessarily $L$-summands, one-sided $L$-summands
are Chebychev, and complete one-sided
    $L$ and $M$-projections are
uniquely determined by their ranges.
These and  other ``one-sided'' analogues of the
classical $M$-ideal theory are presented in \S 3.  We make no attempt to be
exhaustive.  Additional results, together with
a more detailed exposition of the basic theory may be found in
    \cite{Zarikianthesis}.  We have deferred some of these topics to
the sequel of this paper, and to \cite{BESZ}.

Given an operator space $X$ and a completely isometric embedding
\begin{equation}\label{Fembed}
\sigma:X\hookrightarrow B(K,H),
\end{equation}
we say that
$b\in B(H)$ is a \emph{left multiplier} of $X$ if $b\sigma
(X)\subseteq \sigma (X),$
and let
$ M_{\ell }^{\sigma }(X)$ be the algebra of all such $b\in B(H).$ To
simplify the
notation we will often write $X\subseteq B(K,H)$ and $bX\subseteq X$. The left
multipliers in the unital
$C^{*}$-algebra
\[
A_{\ell }^{\sigma }(X)=M_{\ell }^{\sigma }(X)\cap M_{\ell }^{\sigma
}(X)^{*}\subseteq B(H)
\]
are said to be \emph{left adjointable.} Since we have the natural inclusion map
\[
B(K,H)\hookrightarrow B(K\oplus H, K\oplus H)
\]
we may, for most purposes, restrict our attention to multipliers associated
with embeddings of the form $\sigma:X\hookrightarrow B(L)$ for a
Hilbert space
$L$. On the other hand, we need the more general embeddings to prove the
existence of Shilov embeddings (see below).

Given an embedding (\ref{Fembed}), each $b\in M_{\ell }^{\sigma }(X)$
determines a map
\[
\varphi =L^{\sigma }(b):X\rightarrow X:x\mapsto bx,
\]
with $\left\| \varphi \right\| _{cb}\leq \left\| b\right\| .$ We say that a
linear map $\varphi :X\rightarrow X$ is a \emph{left multiplier map} if $
\varphi =L^{\sigma }(b)$ for some embedding $\sigma :X\hookrightarrow B(K,H)$
and $b\in M_{\ell }^{\sigma }(X)\subseteq B(H).$ We let $M_{\ell
}(X)\subseteq CB(X)$ be the set of all such maps $\varphi .$
Similarly, $\varphi $ is
a \emph{left adjointable multiplier map} if $\varphi =L^{\sigma }(b)$
with $b\in
A_{\ell }^{\sigma }(X),$ and we let $A_{\ell }(X)\subseteq CB(X)$
denote the set of
all such maps $\varphi .$

Given an operator space $X$, then one can use the construction of the
``non-commutative Shilov boundary'' of an operator space to find an embedding
$\sigma _{0}:X\hookrightarrow B(K,H)$ with the following properties:
\begin{enumerate}
\item[(i)] for any $\varphi\in M_{\ell }(X)$ there is a unique
element
$b_0\in M_{\ell}^{\sigma_0}(X)$ such that $\varphi=L^{\sigma _{0}}(b_0)$,
\item[(ii)] for any $\varphi\in A_{\ell }(X)$ there exists a unique
element $b_0\in
A_{\ell}^{\sigma_0}(X)$ such that $\varphi=L^{\sigma _{0}}(b_0),$
\item[(iii)] if $\varphi=L^{\sigma}(b_1)$ for some embedding
$\sigma:X\hookrightarrow
B(K_1,H_1)$ and element $b_1\in M_{\ell }^{\sigma}(X)$, then
$\|b_0\|\leq\|b_1\|$.
\end{enumerate}
(see \cite{BlecherShilov}, \cite{Arv1,Arv2}, \cite{Ham1},
\cite{Hamana}). For lack of
a better term, we will refer to an embedding $\sigma_0$ with these properties
as a ``Shilov embedding''. The existence of such an embedding implies
that $M_{\ell
}(X)$ and $A_{\ell }(X)$ are subalgebras of $CB(X).$

If $\sigma_0:X\hookrightarrow B(K,H)$ is a Shilov embedding, then by
definition the
map
\begin{equation}\label{Fmult1}
L^{\sigma _{0}}:M_{\ell }^{\sigma _{0}}(X)\rightarrow M_{\ell }(X)\subseteq
CB(X)
\end{equation}
and its restriction
\begin{equation}
\label{Fmult2}
A_{\ell }^{\sigma _{0}}(X)\rightarrow A_{\ell
}(X)\subseteq CB(X)
\end{equation}
are algebraic isomorphisms.

Since $A_{\ell
}^{\sigma _{0}}(X)$ is a $C^{*}$-algebra, it follows that the algebraic
isomorphism (\ref{Fmult2}) is isometric (see \cite{Tom}, Prop. 1.1),
and we have a
corresponding
$C^{*}$-algebraic structure on $A_{\ell }(X)$. If $\sigma_1$ is another Shilov
embedding, then the algebras $A_{\ell}^{\sigma _{j}}(X)$ $(j=0,1)$
are isometrically
isomorphic as unital Banach algebras. Since a unital norm-decreasing map of
$C^{*}$-algebras is necessarily
$*$-preserving (see, e.g., Lemma \ref{selfadjlem} below), they are isomorphic
$C^{*}\mbox{-algebras}$, and therefore the $C^{*}$-algebraic structure on
$A_{\ell}(X)$ does not depend upon the Shilov embedding. The
self-adjoint projections
in this $C^{*}$-algebra are the complete left $M$-projections on $X$
(see Theorem
\ref{tot2}).

It is shown in \cite{BlecherShilov} and
\cite{BlecherPaulsen} that although the isomorphism
(\ref{Fmult1}) is generally not isometric, there is a natural
operator space structure
on $M_{\ell }(X)$ with respect to which it is an operator algebra. In
particular if
$\varphi\in M_\ell (X)$, then the corresponding norm is given by
$\|\varphi\|_{M_\ell (X)}=\|b_0\|$, where $\varphi=L^{\sigma_0}(b_0)$
for an arbitrary
Shilov embedding $\sigma_0$ and $b_0\in M_\ell^{\sigma_0}(X)$.
 
One of the main objectives of this paper is find
{\em intrinsic} characterizations of the left multiplier and left adjointable
multiplier maps. In order to state  these criteria,
we need some definitions. An element $a$ of a unital Banach
algebra $A$ is said to be \emph{hermitian } if $\left\| e^{ita}\right\| =1$
for all $t\in \Bbb{R}$ (see
\cite{BonsallDuncan}). If $X$ is an operator space,
we say that a mapping $\varphi :X\rightarrow X$ is \emph{completely
hermitian} if it is a hermitian element of $CB(X),$ or equivalently
if, for each
$n\in \Bbb{N},$ the map $\varphi _{n}:M_{n}(X)\rightarrow M_{n}(X)$ is
hermitian in $B(M_{n}(X)).$

We let the space $C_{2}(X)=M_{2,1}(X)$ of $2\times 1$ column matrices over
an operator space $X$ have its canonical operator space structure. Given a
linear mapping $\varphi:X\rightarrow X,$ we define the \emph{column mapping }
$\tau _{\varphi}^{c}:C_{2}(X)\rightarrow C_{2}(X)$ by
\begin{equation*}
\tau _{\varphi }^{c}\left( \left[
\begin{array}{l}
x \\
y
\end{array}
\right] \right) =\left[
\begin{array}{c}
\varphi (x) \\
y
\end{array}
\right] .
\end{equation*}

\begin{theorem}
\label{multiplierth}  Suppose that $X$ is an operator space, and that
    $\varphi:X\rightarrow X$
is a linear mapping. Then the following are equivalent:
\begin{itemize}
\item [(a)]  there  exists
a completely isometric embedding $X\hookrightarrow B(H)$
such that $\varphi (x)=bx$ for some $b\in B(H)$ with
$\left\| b\right\| \leq 1$ (respectively, $b=b^{*}$,
$b$ an orthogonal projection);
     \item [(b)] $\tau_{\varphi }^{c}$ is completely contractive (respectively,
$\tau_{\varphi }^{c}$ is  completely hermitian,
$\varphi $ is a complete left $M$-projection).
\end{itemize}
\end{theorem}

  From our previous discussion of multipliers we may use
a Shilov embedding in (a). It follows that the
first statement in (a) is equivalent to the condition that $\varphi
\in M_\ell(X)$ and
$\|\varphi\|_{M_\ell(X)}\leq 1$.
We will use this result in \S
\ref{applications} to prove that
the left adjointable multiplier algebra
$A_\ell(X)$ of a dual operator space $X$ (i.e. $X$ is the dual of an
operator space) is a von Neumann algebra.
A consequence of this is that
$C^{*}\mbox{-algebraic}$ operator bimodule operations on a dual operator
space are automatically
weak$^{*}$ continuous. We also consider some functorial properties of
the multiplier mappings.

In \S \ref{examples} we give
various examples.  In particular
we prove that the
complete right $M$-ideals in a Hilbert $C^*$-module are exactly
the closed submodules, and we list some consequences of this.  We
also observe that the classical $M$-ideals of Banach spaces, and
the ``complete  $M$-ideals'' of the second author and Ruan, may be viewed
as particular examples of complete left $M$-ideals.

The theory of one-sided ideals and multipliers in a unital $C^{*}$-algebra $
A $ has a long history. It was shown in \cite{Effros} and \cite{Prosser}
    that they are
in one-to-one correspondence with the closed faces of the state space $S(A).$
These faces are particularly well-behaved, and a corresponding theory of
``split faces'' of a convex set was studied in
\cite{AlfsenAndersen,Alfsenbook}. This theory played a key role in the
Alfsen-Schultz characterization of the state spaces of $C^{*}$-algebras
(see \cite{AlfsenSchultz2}).  On
the other hand, K. H. Werner considered a related notion for
\emph{operator systems} (these are matrix \emph{ordered} spaces), and he
defined a notion of multipliers of such spaces \cite{KHWerner}, \cite{Witt}.
E. Kirchberg considered multipliers of a certain class
of operator spaces in
\cite{Kirchberg}.
Arveson was the first to consider ``Shilov
representations''
\cite{Arv1,Arv2}, of operator spaces, and this theory was further developed by
Hamana.  Around 1998, W. Werner considered left
multipliers on a class of {\em non-unital} operator systems and
proved an intrinsic matrix {\em order-theoretic} characterization
which is analogous to our characterization of contractive left
multipliers in Theorem 4.6. Indeed this insightful theorem (in an
early version of \cite{WWerner}) provided the inspiration for our
(non-order theoretic) result. He has very recently pointed out to us
that one can also prove our result by using  a version of the
"Paulsen trick" to replace an operator space by an ordered system of
the variety considered in his paper. In this context one may use a
"$4\times 4$"  matrix argument to recover our theorem. By this trick,
the operator space multipliers in
\cite{BlecherShilov},\cite{BlecherPaulsen} may be described
within Werner's framework, and some of the results from those papers
may be deduced from Werner's work. Similarly the projections that
Werner used in \cite{WWerner} are related to the one-sided
complete $M$-projections of this paper.

\section{Some operator space preliminaries}

\label{prelim}

We refer the reader to the book \cite{EffrosRuanbook} as
    a general reference to the theory of
operator spaces, and for help with any of the details below.

An \emph{operator space} $X$ is a vector space together with distinguished
norms on each matrix space $M_{n}(X)\ $which are linked by the relations
\begin{eqnarray*}
\left\| x\oplus y\right\| &=&\max \left\{ \left\| x\right\| ,\left\|
y\right\| \right\}, \\
\left\| \alpha x\beta \right\| &\leq &\left\| \alpha \right\| \left\|
x\right\| \left\| \beta \right\|.
    \end{eqnarray*}
Here $\alpha, \beta$ are scalar matrices, and the $\oplus$
refers to the ``diagonal direct sum'' of matrices (see \cite{Ruanoper}).
These ``square matrix'' norms
uniquely determine norms on each ``rectangular matrix''
space $M_{m,n}(X).$
By considering matrices over the latter space, we
see that $
M_{m,n}(X)$ is again an operator space. We let
\begin{equation*}
C_{n}(X)=M_{n,1}(X),\,\,R_{n}(X)=M_{1,n}(X),
\end{equation*}
with these operator space structures, and in particular, we let $C_{n}=C_{n}(
\Bbb{C})$ and $R_{n}=R_{n}(\Bbb{C)}.$ We have the natural complete isometries
\begin{eqnarray*}
C_{n}(X) = C_{n}\check{\otimes}X = C_n \otimes_h X, \\
R_{n}(X) = R_{n}\check{\otimes}X = X \otimes_h R_n,  \end{eqnarray*}
where $\check{\otimes}$ and $\otimes_h$ denote the usual
spatial and Haagerup tensor products for operator spaces (see e.g.
chapters 7-9 in
\cite{EffrosRuanbook}). On the other hand, we let
\begin{eqnarray*}
C_{n}[X] = C_{n}\hat{\otimes}X  = X \otimes_h C_n  \\
R_{n}[X] = R_{n}\hat{\otimes}X  = R_n \otimes_h X ,
\end{eqnarray*}
where $\hat{\otimes}$ denotes the projective operator space tensor product.
We have the identifications
\begin{eqnarray*}
(C_{n}(X))^{*} =R_{n}[X^{*}],\,\,\,(R_{n}(X))^{*}=C_{n}[X^{*}]   \; \; ,\\
(C_{n}[X])^{*} =R_{n}(X^{*}),\,\,\,(R_{n}[X])^{*}=C_{n}(X^{*})   \; \; ,
\end{eqnarray*}
where in each case we use the pairings
\begin{equation*}
\left\langle \left[
\begin{array}{l}
x_1 \\
x_2 \\
\vdots \\
    x_n
\end{array}
\right] ,[f_1 \; f_2  \cdots  f_n]\right\rangle =
\sum_k f_k(x_k) =\left\langle [x_1 \; x_2  \cdots  x_n],\left[
\begin{array}{l}
f_1 \\
f_2\\
\vdots \\
    f_n
\end{array}
\right] \right\rangle .
\end{equation*}

An essential distinction between $C_2(X)$ and $C_2[X]$ can be seen from
the following lemma. It should be noted that
the obvious modification of this result is true for rows and
columns of arbitrary length.

\begin{lemma}
Suppose that $X$ is an operator space and that $x,y\in X.$ Then
\begin{equation}
\left\| \left[
\begin{array}{l}
x \\
y
\end{array}
\right] \right\| _{C_{2}(X)}\leq \left( \left\| x\right\| ^{2}+\left\|
y\right\| ^{2}\right) ^{1/2}  \label{ineq1}
\end{equation}
and
\begin{equation}
\left\| \left[
\begin{array}{l}
x \\
y
\end{array}
\right] \right\| _{C_{2}[X]}\geq \left( \left\| x\right\| ^{2}+\left\|
y\right\| ^{2}\right) ^{1/2} \; \; .   \label{ineq2}
\end{equation}
\end{lemma}

\proof
We may assume that $X$ is a subspace of $B(H)$ for some Hilbert space $H.$
Then
\begin{equation*}
\left\| \left[
\begin{array}{l}
x \\
y
\end{array}
\right] \right\| _{C_{2}(X)}^{2}=\left\| \left[
\begin{array}{ll}
x^{*} & y^{*}
\end{array}
\right] \left[
\begin{array}{l}
x \\
y
\end{array}
\right] \right\| =\left\| x^{*}x+y^{*}y\right\| \leq \left\| x\right\|
^{2}+\left\| y\right\| ^{2}.
\end{equation*}
Equivalently, if we let $X\oplus _{2}X$ denote the vector space $X\oplus X$
with the norm
\begin{equation*}
\left\| (x,y)\right\| =(\left\| x\right\| ^{2}+\left\| y\right\| ^{2})^{1/2}
\; ,
     \end{equation*}
then the mapping
\begin{equation*}
\theta _{X}^{c}:X\oplus _{2}X\rightarrow C_{2}(X):(x,y)\mapsto \left[
\begin{array}{l}
x \\
y
\end{array}
\right]
\end{equation*}
is a contraction. Of course the same applies to the corresponding mapping $
\theta _{X}^{r}:X\oplus _{2}X\rightarrow R_{2}(X)$.
If we define
\begin{equation*}
\eta _{X}^{c}:C_{2}[X]\rightarrow X\oplus _{2}X:\left[
\begin{array}{l}
x \\
y
\end{array}
\right] \mapsto (x,y),
\end{equation*}
then it is evident that $(\eta _{X}^{c})^{*}=\theta _{X^{*}}^{r},$ and since
$\theta _{X^{*}}^{r}$ is contractive, that is also true for $\eta _{X}^{c},$
i.e. we have (\ref{ineq2}).
\endproof

It is immediate from the axioms for an operator space that
\begin{equation}
\left\| x\right\| \leq \left\| \left[
\begin{array}{l}
x \\
y
\end{array}
\right] \right\| _{C_{2}(X)},  \label{ineq3}
\end{equation}
and from (\ref{ineq2}) that
\begin{equation}
\left\| x\right\| \leq \left\| \left[
\begin{array}{l}
x \\
y
\end{array}
\right] \right\| _{C_{2}[X]}.  \label{ineq4}
\end{equation}

We will need the following result in Lemma \ref{table1}.

\begin{lemma}  \label{vlem}
Let $X$ be an operator space. Then the map $R_2(R_2[X]) \to
R_2[R_2(X)]$ defined by
\begin{equation}\label{FZarikian}
          \begin{bmatrix} u & v & w & x \end{bmatrix} \mapsto
\begin{bmatrix} u & w & v & x
          \end{bmatrix}
\end{equation}
is a complete isometry.
\end{lemma}
 
\proof
We have from above the natural complete isometries
\begin{eqnarray*}
R_2 \check{\otimes}
(R_2 \hat{\otimes} X) &\cong& (R_2 \hat{\otimes} X) \check{\otimes} R_2\\
                  &\cong& (R_2 \otimes_h X) \otimes_h R_2\\
                  &\cong& R_2 \otimes_h (X \otimes_h R_2)\\
                  &\cong& R_2 \hat{\otimes} (X \check{\otimes} R_2)\\
                  &\cong& R_2 \hat{\otimes} (R_2 \check{\otimes} X).
\end{eqnarray*}
If we successively apply these identifications to an elementary tensor on the
left, we obtain
\begin{eqnarray*}
\begin{bmatrix} \alpha & \beta \end{bmatrix}
\otimes (\begin{bmatrix} \gamma & \delta \end{bmatrix} \otimes x)
    &\mapsto& (\begin{bmatrix} \gamma & \delta \end{bmatrix} \otimes x)
\otimes \begin{bmatrix} \alpha & \beta                    \end{bmatrix}\\
                  &\mapsto& \begin{bmatrix} \gamma & \delta \end{bmatrix}
\otimes (x \otimes \begin{bmatrix} \alpha & \beta
\end{bmatrix})\\
                  &\mapsto& \begin{bmatrix} \gamma & \delta \end{bmatrix}
\otimes (\begin{bmatrix} \alpha & \beta
\end{bmatrix} \otimes x),
\end{eqnarray*}
i.e,
\[
                  \begin{bmatrix} (\alpha \gamma )x & (\alpha\delta)x &
(\beta \gamma )x & (\beta \delta)x
\end{bmatrix} \mapsto
                          \begin{bmatrix} (\gamma \alpha)x & (\gamma
\beta )x & (\delta \alpha)x & (\delta \beta )x
\end{bmatrix}
\]
which coincides with (\ref{FZarikian}). This extends by
linearity to arbitrary tensors on the left.
\endproof

\section{One-sided $M$-projections and $L$-projections}
\label{proj}

If $X$ is a vector space, we say that a linear mapping $P:X\rightarrow X$ is
a \emph{projection }if $P^{2}=P$ (for Hilbert space operators we will also
insist that the mapping be self-adjoint). If $I$ is the identity mapping, it
follows that $I-P$ is also a projection. If $P$ is a projection, then the
linear mappings
\begin{eqnarray*}
\nu_{P}^{c}
&:&X\rightarrow C_{2}(X):x\mapsto
\left[
\begin{array}{c}
P(x) \\
x-P(x)
\end{array}
\right] \; ,  \\
\mu_{P}^{c} &:&
C_{2}(X)\rightarrow X:\left[
\begin{array}{l}
x \\
y
\end{array}
\right] \mapsto P(x)+y-P(y)  \; ,
\end{eqnarray*}
satisfy $\mu _{P}^{c}\circ \nu _{P}^{c}= I$. We have corresponding mappings
$\nu_{P}^{r}:X\rightarrow R_{2}(X)$ and
$\mu_{P}^{r}:R_{2}(X)\rightarrow X$ which satisfy
$\mu_{P}^{r}\circ \nu_{P}^{r}=I$.

We recall that if $X$ is a Banach space, then a projection $P:X\rightarrow X$
is an $M$-\emph{projection} if
for every $x \in X$ we have
\begin{equation*}
\left\| x\right\| =\max \left\{ \left\| P(x)\right\| ,\left\| x-P(x)\right\|
\right\} .
\end{equation*}
If $X$ is an operator space, we say that $P$ is a \emph{complete}
$M$-\emph{projection} if for each $n\in \Bbb{N}$,
$P_{n}:M_{n}(X)\rightarrow M_{n}(X)$
is an $M$-projection.  It
is known that $M$-projections need not be complete $M$-projections
(see \cite{ERcmp}).

    From the introduction, $P:X\rightarrow X$ is a left $M$-projection if and
only if
\begin{equation*}
\nu _{P}^{c}:X\rightarrow C_{2}(X):x\mapsto \left[
\begin{array}{c}
P(x) \\
x-P(x)
\end{array}
\right]
\end{equation*}
is an isometric injection. Using simple matrix manipulations it is evident
that $P$ is a complete left $M$-projection if and only if $\nu _{P}^{c}$ is
completely isometric. Owing to the fact that
\begin{equation*}
\left[
\begin{array}{l}
b \\
a
\end{array}
\right] =\left[
\begin{array}{cc}
0 & 1 \\
1 & 0
\end{array}
\right] \left[
\begin{array}{l}
a \\
b
\end{array}
\right]
\end{equation*}
it is evident that if $P$ is a (complete) left $M$-projection, then the same
is true for $I-P.$

If $e$ is a (self-adjoint) projection in a unital $C^{*}$-algebra $A,$ then $
P(x)=ex$ is a left $M$-projection on $A$ since
\begin{eqnarray}
\left\| \left[
\begin{array}{c}
P(x) \\
x-P(x)
\end{array}
\right] \right\| ^{2} &=&\left\| \left[
\begin{array}{c}
ex \\
x-ex
\end{array}
\right] \right\| ^{2} \label{projection}
=\left\| \left[
\begin{array}{ll}
x^{*}e & x^{*}-x^{*}e
\end{array}
\right] \left[
\begin{array}{c}
ex \\
x-ex
\end{array}
\right] \right\|  \notag \\
&=&\left\| x^{*}ex+x^{*}(1-e)x\right\|  \notag
=\left\| x^{*}x\right\|
=\left\| x\right\| ^{2}. \notag
\end{eqnarray}
If $x\in M_{n}(A),$ then $P_{n}(x)=e_{n}x$ where $e_{n}=e\oplus \cdots
\oplus e$ is a projection in $M_{n}(A),$ and it follows that $P$ is
complete left $M$-projection.

\begin{lemma}\label{2sid}
Suppose that $X$ is an operator space. A projection $P:X\rightarrow X$ is
both a complete left and a complete right $M$-projection if and only if it is a
complete $M$-projection.
\end{lemma}
\proof
If $P$ is a complete left and right $M$-projection, then $\left\| x\right\|
=\left\| \nu _{P}^{c}(x)\right\| =\left\| (\nu _{P}^{r})_{2,1}(\nu
_{P}^{c}(x))\right\| ,$ and thus
\begin{eqnarray*}
\left\| x\right\| &=&\left\| \left[
\begin{array}{c}
P(x) \\
x-P(x)
\end{array}
\right] \right\| =\left\| \left[
\begin{array}{cc}
P^{2}(x) & (I-P)P(x) \\
P(I-P)(x) & (I-P)^{2}(x)
\end{array}
\right] \right\| \\
&=&\left\| \left[
\begin{array}{cc}
P(x) & 0 \\
0 & x-P(x)
\end{array}
\right] \right\| =\max \left\{ \left\| P(x)\right\| ,\left\| x-P(x)\right\|
\right\} .
\end{eqnarray*}
This applies as well to matrices. Conversely if $P$ is a complete
$M$-projection,
then the mapping
\begin{equation*}
\theta :X\rightarrow M_{2}(X):x\mapsto (Px)\oplus (x-Px)
\end{equation*}
is completely isometric.
It follows that
\begin{eqnarray*}
\left\| \left[
\begin{array}{c}
P(x) \\
x-P(x)
\end{array}
\right] \right\| &=&\left\| \theta _{2,1}\left( \left[
\begin{array}{c}
P(x) \\
x-P(x)
\end{array}
\right] \right) \right\| \\
&=&\max \left\{ \left\| Px\right\| ,\left\| x-Px\right\| \right\} =\left\|
x\right\| .
\end{eqnarray*}
Again it is easy to generalize this to matrices. A similar argument may be
applied to row matrices.
\endproof

The following result will be useful in our discussion of duality.

\begin{proposition}
\label{munu}If $X$ is an operator space and $P:X\rightarrow X$ is a
projection, then $P$ is a complete left $M$-projection if and only if $\mu
_{P}^{c}$ and $\nu _{P}^{c}$ are both completely contractive.
\end{proposition}
\proof
If $\nu _{P}^{c}$ is completely isometric, then
\begin{equation*}
\left\| P(x)+y-P(y)\right\| =\ \left\| \left[
\begin{array}{c}
P(x) \\
y-P(y)
\end{array}
\right] \right\| \leq \left\| \left[
\begin{array}{c}
P(x) \\
x-P(x) \\
P(y) \\
y-P(y)
\end{array}
\right] \right\| =\left\| \left[
\begin{array}{l}
x \\
y
\end{array}
\right] \right\| ,
\end{equation*}
and thus $\mu _{P}^{c}$ is
contractive. These calculations work
as well for matrices. The converse is trivial
since if two complete contractions compose to the identity, then
the first is completely isometric.
\endproof

As in the Banach space theory, left $M$-projections have certain automatic
continuity properties.

\begin{proposition}
\label{weakstarcont}Suppose that $X$ is
an  operator space which is also a dual Banach space.
Then any left $M$-projection $P:X\rightarrow X$ is weak$^{*}$
continuous.
\end{proposition}

\proof
A standard argument in functional analysis
shows that it
suffices to prove that the unit balls of $P(X)$ and
$(I-P)(X)$ are weak$^{*}$ closed.  By symmetry of $P$ and $I-P$
it is enough to prove the former.
Let us suppose that $\{ y_{\nu } \}$ is a
net in
$P(X)$ with $\left\| y_{\nu }\right\| \leq 1,$
converging weak$^{*}$ to an element
$x \in X.$ If we let $y=P(x)$ and $z=(I-P)(x),$ it follows that $y_{\nu
}^{\prime }=y_{\nu }-y$ converges weak$^{*}$ to
$z.$  Scaling by $\frac{1}{2}$, we may suppose that we have
a net $\left\| y_{\nu }\right\| \leq 1,$ converging weak$^{*}$ to
a  $z \in (I-P)(X)$.
For any $t>0,$ we have $y_\nu + t z \rightarrow (1+t) z$.  Hence
using the fact that norm closed balls in $X$ are weak$^{*}$ closed,
and  (\ref{ineq1}), we see that $$ (1+t)^2 \Vert z \Vert^2
    \; \leq \; \sup_\nu \;
\Vert y_\nu + t z \Vert^2 \;
    = \;  \sup_\nu \; \left\| \left[
\begin{array}{l} y_\nu \\ tz \end{array}
\right] \right\|^2_{C_2(X)} \;
    \leq \; 1 + t^2 \Vert z \Vert^2  \; \; . $$
Letting $t \rightarrow \infty$ shows that $z = 0$.\endproof

We say that a projection $P:X\rightarrow X$ is a \emph{left}
$L$\emph{-projection} if
\begin{equation*}
\left\| x\right\| =\left\| \left[
\begin{array}{c}
P(x) \\
x-P(x)
\end{array}
\right] \right\| _{C_{2}[X]},
\end{equation*}
or equivalently
if $\nu _{P}^{c}:X\rightarrow C_{2}[X]$ is isometric. We say
that $P$ is a \emph{complete left} $L$-\emph{projection} if the mapping $\nu
_{P}^{c}:X\rightarrow C_{2}[X]$ is a complete isometry.

\begin{proposition}
\label{Lproj}If $X$ is an operator space and $P:X\rightarrow X$ is a
projection, then $P$ is a complete left $L$-projection if and only if $\nu
_{P}^{c}:X\rightarrow C_{2}[X]$ and $\mu _{P}^{c}:C_{2}[X]$ $\rightarrow X$
are completely contractive.
\end{proposition}
\proof
As in the proof of Proposition 3.2, the key point is to show that if $P$ is a
complete left $L$-projection, i.e. $\nu _{P}^{c}$ is completely isometric,
then $\mu_{P}^{c}$ is a complete contraction. The truncation mapping
\[
\rho :C_{4}\rightarrow C_{2}:\left[
\begin{array}{l}
\alpha  \\
\beta  \\
\gamma  \\
\delta
\end{array}
\right] \mapsto \left[
\begin{array}{l}
\alpha  \\
\delta
\end{array}
\right]
\]
is completely contractive and thus, by the
``functoriality'' of the projective tensor
product, it induces a
complete contraction $\rho
\otimes id:C_{4}[X]\rightarrow C_{2}[X].$ We have a
commutative diagram
\[
\begin{CD}
C_{4}[X]            @>\rho \otimes id>>      C_{2}[X]     \\
         @AA\nu'_{P}A              @AA\nu _{P}^{c}A \\
       C_{2}[X]          @>\mu _{P}^{c}>>       X
\end{CD}
\]
where $\nu'_P = id_{C_{2}}\otimes \nu^c_{P}$.  This
is because for any $x,y\in X,$
\[
(\rho \otimes id)\circ \nu'_{P}\left[
\begin{array}{l}
x \\
y
\end{array}
\right] = (\rho  \otimes id) \left[ \begin{array}{l}
Px \\
x-Px \\
Py \\
y-Py
\end{array}
\right] =\left[
\begin{array}{l}
Px \\
y-Py
\end{array}
\right] =\nu^c_{P}(Px+y-Py).
\]
It follows that $(\rho \otimes id)\circ \nu'_{P}$ has
range in $\nu^c_{P}(X).$ By hypothesis,
\[
\nu^c_{P}:X\rightarrow \nu^c_{P}(X)
\]
is a complete isometry, and thus
\[
\mu^c_{P}= (\nu^c_{P})^{-1}\circ (\rho \otimes id) \circ \nu'_{P}
\]
is a complete contraction.
\endproof

The proofs of Proposition \ref{munu} and
Proposition \ref{Lproj} do not generalize
to left $M$- and left $L$-projections.
For this reason it might be useful to consider a related
notion. We say that a
projection  $P:X\rightarrow X$ is a \emph{strong left}
$M$\emph{-projection} if $\nu _{P}^{c} : X \rightarrow C_2(X)$
and $\mu _{P}^{c} :  C_2(X) \rightarrow  X$ are contractive,
and we similarly define \emph{strong left} $L$\emph{-projections}.  The
reader will see that
the duality relationships considered below are also valid for
these ``strong'' one-sided projections.  In fact
most of the results of this section
which are stated for ``complete one-sided projections and summands
and ideals'', are also valid
with ``complete'' replaced by ``strong''.
 
\begin{corollary}
\label{LMduality}If $X$ is an operator space and $P:X\rightarrow X$ is a
projection, then $P$ is a
complete left $M$-projection if and only if $P^{*}$ is a
complete right $L$-projection. Similarly $P$ is a
complete right $L$-projection if and only if $P^{*}$ is a
complete left $M$-projection.
\end{corollary}

\proof
For any $x\in X$ and $f,g\in X^{*},$
\begin{eqnarray*}
\langle \nu _{P}^{c}(x),[f\,\,g]\rangle &=&\left\langle \left[
\begin{array}{c}
P(x) \\
x-P(x)
\end{array}
\right] ,[f\,\,g]\right\rangle \\
&=&\langle P(x),f\rangle +\langle x-P(x),g\rangle \\
&=&\langle x,P^{*}(f)+g-P^{*}(g)\rangle \\
&=&\langle x,\mu _{P^{*}}^{r}([f \; g]) \rangle \; ,
\end{eqnarray*}
and thus $(\nu _{P}^{c})^{*}=\mu _{P^{*}}^{r}.$ Similarly, $(\mu
_{P}^{c})^{*}=\nu _{P^{*}}^{r}.$ It follows from Proposition \ref{munu} and
Proposition \ref{Lproj}, and basic operator space duality,
    that $P$ is a complete left $M$-projection if and only if $P^{*}$
is a complete right $L$-projection, and similarly $P$ is a complete right
    $L$-projection if and only if $P^{*}$ is a complete left $M$-projection.
\endproof

We recall from the
introduction that a subspace $J$ of an operator space $X$
    is a\emph{\ (complete) right} $M$-\emph{summand} of $X$ if it is the
range
of a (complete) left $M$-projection.
We say that $J$ is a \emph{(complete)} \emph{right} $L$\emph{-summand}
if it is the range of a (complete) left $L$-projection$.$ We
note that if $P:X\rightarrow X$ is a bounded projection, then the same is
true for $P^{*},$ and we have
\begin{equation*}
P(X)^{\perp }=\ker P^{*}=(I-P^*)(X^{*}).
\end{equation*}
We thus have

\begin{corollary}
If $X$ is an operator space and $J\subseteq X$ is a complete right
$M\mbox{-sum}$\-mand, then $J^{\perp }$ is a complete left
$L$-summand, and if $
J\subseteq X$ is a complete right $L$-summand, then $J^{\perp }$ is a
complete left $M$-summand.
\end{corollary}

A subspace $J$ of a Banach space $X\ $is said to be \emph{proximinal}
(respectively, \emph{Chebychev)} if for each $x\in X,$ the set
\begin{equation*}
\mathcal{P}_{J}(x)=\left\{ h\in
J:\left\| x-h\right\| =\left\| x-J\right\|
\right\}
\end{equation*}
is non-empty (respectively, has one point). If $P$:$X\rightarrow X$ is a
left $M$-projection, $J=P(X),$ and $x\in X,$ then
\begin{equation*}
P(x)\in \mathcal{P}_{J}(x),
\end{equation*}
since if $x\in X$ and $h\in J,$ then
\begin{equation*}
\left\| x-h\right\| =\left\| \left[
\begin{array}{c}
P(x-h) \\
(I-P)(x-h)
\end{array}
\right] \right\| =\left\| \left[
\begin{array}{l}
P(x)-h \\
x-P(x)
\end{array}
\right] \right\| \geq \left\| x-P(x)\right\| .
\end{equation*}
It follows that right $M$-summands are proximinal. A similar argument
with (\ref{ineq4}) shows that right $L$-summands
are also proximinal.

\begin{proposition}
If $P$ is a left $L$-projection with $J=P(X)$, then $\mathcal{P}
_{J}(x)=\{P(x)\},$ and thus $J$ is Chebychev.
\end{proposition}
\proof
If $h\in J$, then from (\ref{ineq2}),
\begin{eqnarray*}
\left\| x-h\right\| ^{2} &=&\left\| \left[
\begin{array}{c}
P(x-h) \\
\,(I-P)(x-h)
\end{array}
\right] \,\,\right\| _{C_{2}[X]}^{2} \\
&=&\left\| \left[
\begin{array}{l}
P(x)-h \\
\,x-P(x)
\end{array}
\right] \,\,\,\right\| _{C_{2}[X]}^{2} \\
&\geq &\left\| P(x)-h\right\| ^{2}+\left\| x-P(x)\right\| ^{2}.
\end{eqnarray*}
It follows that if $h\in \mathcal{P}_{J}(x)$ $,$ then $\left\| x-h\right\|
=\left\| x-P(x)\right\| $ and $h=P(x).$
\endproof

\begin{corollary}
If $J$ is a
complete right $M$-summand (respectively, right
$L\mbox{-sum}$\-mand), then there
is only one complete
left $M$-projection (respectively, left $L\mbox{-pro}$\-jection) with
range $J.$
\end{corollary}

\proof Given
left $L$-projections $P$ and $Q$ with
$J=P(X)=Q(X),$ we have
\begin{equation*}
\left\{ P(x)\right\} =\mathcal{P}_{J}(x)=\left\{ Q(x)\right\}
\end{equation*}
for $x \in X$,
and therefore $P=Q.$ If $P$ and $Q$ are complete
left $M$-projections with $
J=P(X)=Q(X),$ then
\begin{equation*}
\ker P^{*}=J^{\perp }=\ker Q^{*}
\end{equation*}
implies that the right $L$-projections $I-P^{*}$ and $I-Q^{*}$ have the same
range. Thus $I-P^{*}=I-Q^{*}$ and $P=Q.$
\endproof

In the introduction we defined a subspace $J$ of an operator space $X$ to be
a right $M$-ideal if $J^{\perp \perp }$ is a right $M$-summand. From
the next result we see that
it is equivalent to assume that $J^{\perp }$ is a left
$L$-summand. As in
    the Banach space theory,
this next result also shows that there is no need to define
$L$-ideals, since they must coincide with $L$-summands.
 
\begin{proposition}
If $J$ is a closed subspace of an operator space $X$ for which $J^{\perp }$
is a
complete right $M$-summand, then $J$ is a complete
left $L$-summand.   Indeed
any complete right $M$-summand in a dual operator space $X^*$
is the annihilator of a complete left $L$-summand in $X$.
     \end{proposition}
 
\proof
Let us suppose that $J^{\perp }$ is a
complete right $M$-summand in $X^{*}\ $and let $P
$ be the complete
left $M$-projection onto $J^{\perp }.$ From Proposition \ref
{weakstarcont}, $P$ is weak$^{*}$ continuous. It follows that $P=Q^{*}$ for a
projection $Q:X\rightarrow X.$ That implies that
\begin{equation*}
(I-Q)(X)^{\perp }=\ker (I-P)=P(X^{*})=J^{\perp }
\end{equation*}
and thus $J=(I-Q)(X).$ Since $Q^{*}$ is a
complete left $M$-projection, $Q$ and $I-Q$
are
complete right  $L$-projections.  The proof for the second assertion
is similar.
\endproof

In fact stronger versions of the last few results are true.  We omit
the proofs, which are very simple and identical to their classical
versions (see \cite{HWW}):

\begin{theorem} \label{list}   In the following, $X$ is an operator space.
\begin{itemize}
\item[(a)]   Suppose that $P$ is a complete left
$M$-projection on $X$.
If $Q$ is a contractive projection on $X$
with $Ran \; Q = Ran \; P$, then $P = Q$.
\item[(b)]  Suppose that $P$ is  a complete right $L$-projection
    on  $X$.  If  $Q$ is a
contractive projection on $X$ with
    $ker \; Q = ker \; P$, then $Q = P$.
\item[(c)]   If there exists a contractive projection
from $X$
onto a complete
right $M$-ideal $J$ of  $X$, then $J$ is  a complete
right $M$-summand.
Moreover such a contractive projection is then unique.
\item[(d)]   If $J$ is a complete
    right $M$-ideal in $X$, and if
$J$ is a dual Banach space, then $J$ is a
complete right $M$-summand in $X$.
\item[(e)]  If $X$ is a dual operator space, and if $J$
is a weak*-closed complete
right $M$-ideal of $X$, then $J$ is a complete right $M$-summand in $X$ which
is the annihilator of a complete left $L$-summand in $X_*$.
\end{itemize}
\end{theorem}

The ``complete'' hypothesis in the  results above may be weakened
to the ``strong'' condition briefly alluded to earlier. In light of
the topics to
be discussed in \S 6.5 below, (e) may be regarded as an operator
space generalization of the
result that  weak*-closed submodules of self-dual $C^*$-modules are
orthogonally
complemented. (d) is related to the well-known fact
that if a closed submodule of a Hilbert $C^{*}\mbox{-mod}$\-ule is
self-dual, then it is orthogonally complemented.

\begin{proposition}
Suppose that $X$ and $Y$ are operator spaces. If $P$ is a
complete left $M$-projection on $X,$ then
\begin{equation*}
P\otimes id:X\check{\otimes}Y\rightarrow X\check{\otimes}Y
\end{equation*}
is a complete
left $M$-projection. If $P$ is a
complete left $L$-projection on $X$,
then
\begin{equation*}
P\otimes id:X\hat{\otimes}Y\rightarrow X\hat{\otimes}Y
\end{equation*}
is a complete left $L$-projection.
\end{proposition}
\proof
Owing to the functorial properties of the tensor product,
and using Proposition \ref{munu},
the mappings $\mu^c_P$ and $\nu^c_P$ tensor with $id_Y$ to give
complete contractions
\begin{equation*}
X\check{\otimes}Y\rightarrow C_{2}(X)\check{\otimes}Y\rightarrow X\check{
\otimes}Y.
\end{equation*}
The first relation then follows
again from Proposition \ref{munu}, together with the simple identification
\begin{equation*}
C_{2}(X\check{\otimes}Y)=C_{2}(X)\check{\otimes}Y.
\end{equation*}
The second relation
follows similarly.
\endproof

\section{Multipliers}

\label{mult}

In order to illustrate the definition of the multiplier mappings, let us
consider an elementary proof for
the characterization of complete left $M$-projections given in
Theorem \ref{multiplierth}.

\begin{proposition}
\label{Mproject}The complete left $M$-projections in an operator space $X$
are just the mappings $P(x) = ex$ for a completely isometric
embedding $X \hookrightarrow B(H)$ and an orthogonal
projection $e \in B(H)$.
\end{proposition}

\proof
If $P:X\rightarrow X$ is a complete left $M$-projection, then let us
fix an embedding
$X\subseteq B(H).$ By definition, the mapping
\begin{equation*}
\sigma :X\hookrightarrow B(H\oplus H):x\mapsto \left[
\begin{array}{cc}
P(x) & 0 \\
(I-P)(x) & 0
\end{array}
\right]
\end{equation*}
is completely isometric. We have that
\begin{equation*}
\sigma (P(x))=\left[
\begin{array}{cc}
P(x) & 0 \\
0 & 0
\end{array}
\right] =\left[
\begin{array}{cc}
1 & 0 \\
0 & 0
\end{array}
\right] \sigma(x),
\end{equation*}
and thus $e=\left[
\begin{array}{cc}
1 & 0 \\
0 & 0
\end{array}
\right] \in B(H\oplus H)$ is the desired left projection relative to the
embedding $\sigma $.
The converse is immediate (see the calculation before Lemma
\ref{2sid}).
\endproof
We will give some other characterizations
of the complete left $M$-projections in Theorem \ref{tot2}.

In order to prove the remaining parts of
Theorem \ref{multiplierth}, it is useful to consider a
bimodule version of Hamana's theory of injective envelopes
\cite{BlecherPaulsen}.
Given unital $C^{*}\mbox{-al}$ge\-bras $A$ and $B,$ an operator
space $X$
which is also a left $A$-module  is called a \emph{left operator
    }$A$\emph{-module} if $\left\|
ax\right\| \leq \left\| a\right\| \left\| x\right\| $ for all matrices $a\in
M_{n}(A)$ and $x\in M_{n}(X).$
We assume that the module action is {\em unitary}, i.e. that
$1x=x$ for all $x$. There is a similar definition for right operator
$B$-modules, and for
operator $A$--$B$-modules.
Bimodule mappings are defined in the usual manner.

We say that an
operator  $A$--$B$-bimodule $Z$ is an \emph{injective}
bimodule if given an inclusion
of $A$--$B$-bimodules $X\subseteq Y$, any
completely contractive $A$--$B$-bimodule mapping
$\theta :X\rightarrow Z$
extends to
an $A$--$B$-bimodule mapping
$Y\rightarrow Z$. An inclusion of $A$--$B$-bimodules
$X\subseteq Y$ is \emph{rigid }if given a completely
contractive $A$--$B$-bimodule mapping $\varphi :Y\rightarrow
Y$ such that $\varphi _{|X}=id_{X},$ it follows that $\varphi
=id_{Y}.$
We say that an injective operator $A$--$B$-bimodule $Z$ is an
    operator
$A$--$B$-{\em bimodule injective envelope} of an operator
$A$--$B$-bimodule $X$,   if there exists a completely isometric
rigid $A$--$B$-bimodule inclusion $X \hookrightarrow
Z$.  Following Hamana's argument \cite{Ham1,Hamana},
one can see that  the $A$--$B$-bimodule injective envelope is unique
in the obvious sense.  If $A=B=\Bbb{C},$ then
we are simply talking about the
injective envelope
$I(X)$ of an operator space $X$, as discussed in
\cite{Ham1,Hamana,Ruaninjec}.
The following result was proved in
\cite{BlecherPaulsen}, Corollary
2.6.  In fact we only need the
${\Bbb C}\oplus {\Bbb C}$--$\,{\Bbb C}$-module version
of this result, which may be proved by elementary methods.

\begin{lemma}
\label{injecmod}  An
operator $A$--$B$-bimodule $Y$ is injective as an
operator $A$--$B\mbox{-bimod}$\-ule if and only if it is injective as
an operator space. The
injective envelope $I(X)$ of the operator space $X$ may be regarded as the
operator $A$--$B$-bimodule injective envelope of $X.$
\end{lemma}

We will
be considering infinite matrices over operator spaces. Given
an operator space $X\ $and cardinals $m,n,$ we have a corresponding operator
space $M_{m,n}(X)$ of all matrices for which the finite truncations are
uniformly bounded (see \cite{EffrosRuanbook}).
If $\varphi :X\rightarrow Y$ is a
completely bounded mapping of operator spaces, the mapping
\begin{equation*}
\varphi _{m,n}:M_{m,n}(X)\rightarrow M_{m,n}(Y):\left[ x_{ij}\right]
\rightarrow \left[ \varphi (x_{ij})\right]
\end{equation*}
satisfies
$\left\| \varphi_{m,n} \right\| _{cb}=\left\| \varphi \right\|
_{cb}.$
If we let $D_{n}$ denote the diagonal matrices in $M_{n}$ it is
evident that $M_{m,n}(X)$ is an operator
   $D_{m}$--$D_{n}$-bimodule. The $D_{m}$--$D_{n}$-bimodule mappings
\begin{equation*}
\varphi :M_{m,n}(X)\rightarrow M_{m,n}(Y),
\end{equation*}
are just those for which there exist linear mappings $\varphi
_{ij}:X\rightarrow Y$ with
\begin{equation*}
\varphi ([x_{ij}])=[\varphi _{ij}(x_{ij})].
\end{equation*}

We will only need the following result for $m = 2, n = 1,$ in which
case there is also an elementary direct proof.  We have included the general
case since it is of independent interest.
\begin{lemma}
\label{matrixinj}For any cardinals $m,n,$ we have a natural identification
\begin{equation*}
M_{m,n}(I(X))\cong I(M_{m,n}(X)),
\end{equation*}

i.e. $M_{m,n}(I(X))$ is an injective envelope of
$M_{m,n}(X)$.  \end{lemma}

\proof From
the previous lemma
it suffices
to prove that $M_{m,n}(I(X))$ is the
$D_{m}$--$D_{n}$-bimodule
injective envelope of $M_{m,n}(X)$. To see this we
first note that if $Z$ is injective, then so
is $M_{m,n}(Z)$.
This follows since if $\pi :B(H)\rightarrow Z$ is a surjective
completely  contractive projection, then
\begin{equation*}
\pi _{m,n}:M_{m,n}(B(H))\rightarrow M_{m,n}(Z)
    \end{equation*}
is a completely contractive projection of the injective operator
space $M_{m,n}(B(H))
\cong B(H^n,H^m)$ onto $M_{m,n}(Z)$. If $\varphi :M_{m,n}(I(X))\rightarrow
M_{m,n}(I(X))$ is a
$D_{m}$--$D_{n}$-bimodule complete
    contraction such that $\varphi _{|M_{m,n}(X)}=id_{M_{m,n}(X)},$ then in
particular, $\varphi_{ij}(x)=x$ for $x\in X,$ and therefore $\varphi
_{ij}(x)=x $ for $x\in I(X).$ It follows that $\varphi =id$ and we see that
\begin{equation*}
M_{m,n}(X)\subseteq M_{m,n}(I(X))
\end{equation*}
is a rigid bimodule inclusion. Thus $M_{m,n}(I(X))$ is a bimodule
injective envelope of $M_{m,n}(X)$ and from Lemma \ref{injecmod}
it is an operator space injective
envelope of $M_{m,n}(X).$
\endproof

\begin{lemma}
\label{bigmorita}Suppose that $X$
is the second dual of a ternary system.  Then for some
cardinal $J,$ $M_{J}(X)$ is completely isometric to a
von Neumann algebra.
\end{lemma}

    This result is in the
    folklore of the Morita
equivalence theory of von Neumann algebras. It may be
found in \cite{BlecherShilov} Lemma 5.8, and a more general result
assuming that  $X$  is a weakly closed injective ternary system
may also be deduced from results in \cite{EOR}.

We will use the following simple but elegant result of R. R. Smith (see
\cite{BESZ}).   We include a sketch of
the proof for the sake of completeness.

\begin{lemma}
\label{vonNeum} \cite{BESZ}
    Suppose that $M$ is a von Neumann algebra. Then a mapping $
\varphi :M\rightarrow M$ has the form $\varphi (x)=bx$ for some $b\in M$
with $\left\| b\right\| \leq 1$ if and only if the column mapping
\begin{equation*}
\tau _{\varphi }^{c}:C_{2}(M)\rightarrow C_{2}(M):\left[
\begin{array}{l}
x \\
y
\end{array}
\right] \mapsto \left[
\begin{array}{c}
\varphi (x) \\
y
\end{array}
\right]
\end{equation*}
is contractive.
\end{lemma}

\proof
For the difficult direction, we
suppose that
$\tau _{\varphi }^{c}$ is contractive, and apply
$\tau _{\varphi }^{c}$ to the column in $C_2(M)$ with entries
$e$ and $1-e$, for an orthogonal
projection $e \in M$.   We obtain
$\varphi (e)^{*}\varphi (e)+(1-e)\leq 1$ and thus
\begin{equation*}
(1-e)\varphi (e)^{*}\varphi (e)(1-e)=0,
\end{equation*} giving
$\varphi (e)(1-e)=0.$ But this
relation also holds for the projection $1-e,$
i.e., we have $\varphi (1-e)e=0.$ We conclude that
\begin{equation*}
\varphi (e)=\varphi (e)e=\varphi (1)e.
\end{equation*}
Since the linear span of the
projections is norm dense in $M,$ $\varphi
(x)=bx$ for all $x\in M$,
where $b = \varphi (1)$.
\endproof

\begin{theorem}
\label{lmu}
Suppose that $X$ is an operator space and that $\varphi :X\rightarrow X$ is
a linear mapping. Then there exists a completely isometric
    embedding $X \hookrightarrow B(H)$ and an operator
$b\in B(H)_1$ with $\varphi (x)=bx$ for all $x \in X$ if and only if
\begin{equation}
\tau _{\varphi }^{c}:C_{2}(X)\rightarrow C_{2}(X):\left[
\begin{array}{l}
x \\
y
\end{array}
\right] \mapsto \left[
\begin{array}{c}
\varphi (x) \\
y
\end{array}
\right]   \label{columncontr}
\end{equation}
is completely contractive.
\end{theorem}

\proof
Let us suppose that $\tau _{\varphi }^{c}$ is completely contractive. From
Lemma \ref{matrixinj} and Lemma \ref{injecmod} , $C_{2}(I(X))=I(C_{2}(X))$
is the
$D_{2}$--$D_{1}$-bimodule injective envelope of $C_{2}(X).$ Thus we
may extend the
$D_{2}$--$D_{1}$-bimodule mapping
\begin{equation*}
\tau _{\varphi }^{c}:C_{2}(X)\rightarrow C_{2}(X)
\end{equation*}
to a bimodule mapping
\begin{equation*}
\theta :C_{2}(I(X))\rightarrow C_{2}(I(X)):\left[
\begin{array}{l}
x \\
y
\end{array}
\right] \mapsto \left[
\begin{array}{c}
\theta _{1}(x) \\
\theta _{2}(y)
\end{array}
\right] .
\end{equation*}
Since $\theta _{2}$ restricts to the identity on $X$ and $X\subseteq I(X)$
is rigid, $\theta _{2}=id_{I(X)}.$ Thus if we let $\tilde{\varphi}=\theta
_{1}:I(X)\rightarrow I(X),$ it follows that
\begin{equation*}
\tau _{\tilde{\varphi}}^{c}:C_{2}(I(X))\rightarrow C_{2}(I(X)):\left[
\begin{array}{l}
x \\
y
\end{array}
\right] \mapsto \left[
\begin{array}{c}
\tilde{\varphi}(x) \\
y
\end{array}
\right]
\end{equation*}
is completely contractive.

If we use the natural identification
\begin{equation*}
C_{2}(I(X)^{**})=C_{2}(I(X))^{**},
\end{equation*}
it follows that
\begin{equation*}
\tau _{\tilde{\varphi}^{**}}^{c}=(\tau _{\tilde{\varphi}
}^{c})^{**}:C_{2}(I(X)^{**})\rightarrow C_{2}(I(X)^{**}):\left[
\begin{array}{l}
x \\
y
\end{array}
\right] \mapsto \left[
\begin{array}{c}
\tilde{\varphi}^{**}(x) \\
y
\end{array}
\right]
\end{equation*}
is completely contractive. We have
by \cite{Hamana,Ruaninjec} that $I(X)$ is completely isometric to a
ternary system $eA(1-e),$ where $A$ is a $C^{*}$-algebra and $e$ is
an orthogonal
projection in $A.$ It follows that $I(X)^{**}$ is completely
isometric to the weakly
closed ternary system $eA^{**}(1-e),$ and from Lemma \ref{bigmorita} there
is a cardinal $J$ such that $R=M_{J}(I(X)^{**})$ is a
von Neumann algebra.
The corresponding mapping
\begin{equation*}
\bar{\varphi}=(\tilde{\varphi}^{**})_{J}:M_{J}(I(X)^{**})\rightarrow
M_{J}(I(X)^{**})
\end{equation*}
extends the mapping
\begin{equation*}
\varphi _{J}:M_{J}(X)\rightarrow M_{J}(X),
\end{equation*}
and from the identification
\begin{equation*}
C_{2}(M_{J}(I(X)^{**}))=M_{J}(C_{2}(I(X)^{**}))
\end{equation*}
we have that
\begin{equation*}
\tau _{\bar{\varphi}}^{c}:C_{2}(R)\rightarrow
C_{2}(R):\left[
\begin{array}{l}
x \\
y
\end{array}
\right] \mapsto \left[
\begin{array}{c}
\bar{\varphi}(x) \\
y
\end{array}
\right]
\end{equation*}
is completely contractive.

    From Lemma \ref{vonNeum}, we have that there is a contraction
$b\in R$ such
that $\bar{\varphi}(x)=bx$ for all $x \in R$.
Let us fix an index $j_{0}\in J,$ and if $x\in X,
$ define
\begin{equation*}
\lbrack x]_{j_{0}}\in M_{J}(X)\subseteq M_{J}(I(X)^{**})
\end{equation*}
to be the matrix with $x$ at the $j_{0},j_{0}$ entry and zero elsewhere.
Then
\begin{equation*}
\lbrack \varphi (x)]_{j_{0}}=\varphi _{J}([x]_{j_{0}})=\bar{\varphi}
([x]_{j_{0}})=b[x]_{j_{0}}.
\end{equation*}
The last product here needs a word of clarification.  The point is that
$[x]_{j_{0}}$ is in $M_J(I(X)^{**})$ which is only linearly
completely isometric, via a mapping $\rho$ say,
    to the von Neumann algebra $R$.  Then the statement above
reads, more precisely,
$$\rho(\lbrack \varphi (x)]_{j_{0}}) = b \rho([x]_{j_{0}}) \; \; . $$
Defining an
embedding of $X$ in $R$ by  $\sigma_{1}(x) = \rho([x]_{j_{0}})$, we
see that $\varphi$ is a left multiplier mapping.

We leave the simple argument for the converse to the reader.
\endproof
 
We remark that the mapping $\sigma_{1}$ constructed in the previous
proof cannot
take the place of the Shilov embedding $\sigma_{0}$ described in the
introduction, since in
particular the  corresponding mapping $L : M_{\ell}^{\sigma_{1}}(X)
\rightarrow CB(X)$
is not one-to-one. On the other hand with a little effort, and using
results in \cite{BlecherPaulsen}, it may be seen that a compression of
$\sigma_{1}$ has the desired properties of $\sigma_{0}$. The space of
relative multipliers with respect to this compression
will then coincide with the $IM_{\ell}(X)$ formulation of
the left multiplier algebra given in \cite{BlecherPaulsen}.

The procedure used in the
last proof of passing from $X$ to $I(X)$ to $I(X)^{**}$ and finally to
the von Neumann algebra $R \cong M_J(I(X)^{**})$ was first used in
\cite{BlecherShilov} \S 5.
These steps provide a useful and essentially canonical technique for
embedding an arbitrary operator space
$X$ into a von Neumann algebra.

\begin{corollary} \label{tot}  If $\varphi$ is a linear
mapping on a right $C^*$-module, then $\varphi$ is a
contractive module mapping if and only if $\tau^c_{\varphi}$ is
completely contractive.   \end{corollary}

\proof
This  follows from  the last theorem
and the fact from \cite{BlecherShilov}
A.4 that $M_\ell(Z)$ for a
$C^*$-module $Z$ is the set of bounded module mappings on $Z$.
\endproof

\begin{corollary}
\label{unit}
Suppose that $X$ is an operator space and that $\varphi :X\rightarrow X$ is
a linear mapping. Then $\varphi (x)=
ux$ for a unitary $u\in
A_{\ell }(X)$ if
and only if $\tau _{\varphi }^{c}$ is a completely isometric bijection.
\end{corollary}

\proof
One direction is clear.  For the other, let us use a
Shilov embedding
$\sigma_0 : X \hookrightarrow B(K,H)$
(see \S \ref{Int}). Then applying Theorem \ref{lmu} to
$\varphi $ and $\varphi^{-1},$ we obtain contractions $b,c\in
B(H)$ with $bcx=x=cbx.$
Since $L^{\sigma_0} : M_\ell^{\sigma_0}(X)
\rightarrow CB(X)$ is one-to-one, it follows that $
bc=cb=1,$ and thus $b=c^{-1}.$ Since $b$ and $c$ are both contractions, $b$
is unitary.
\endproof

We can now prove the remaining  assertion  in Theorem \ref{multiplierth},
namely the characterization of left self-adjointable multipliers.

\begin{theorem}
\label{cch}
If $X$ is an operator space, then a mapping $\varphi :X\rightarrow X$ is a
left self-adjoint multiplier if and only if
$\tau^c_\varphi$ is completely hermitian.
\end{theorem}

\proof
One direction is fairly clear.  For the other, we have that
\begin{eqnarray*}
\exp it\tau _{\varphi }^{c} \left( \left[
\begin{array}{l}
x \\
y
\end{array}
\right] \right)
    &=&(I+it\tau _{\varphi }^{c}+\frac{(it\tau _{\varphi }^{c})^{2}}{2!}
+\cdots )\left( \left[
\begin{array}{l}
x \\
y
\end{array}
\right] \right)
\\
&=&\left[
\begin{array}{c}
\exp it\varphi (x) \\
e^{it}y
\end{array}
\right] \\
&=&\left[
\begin{array}{cc}
1 & 0 \\
0 & e^{it}
\end{array}
\right] \left[
\begin{array}{c}
\exp it\varphi (x) \\
y
\end{array}
\right]
\end{eqnarray*}
and thus
\begin{equation*}
\exp it\tau _{\varphi }^{c}=\left[
\begin{array}{cc}
1 & 0 \\
0 & e^{it}
\end{array}
\right] \tau _{\exp it\varphi }^{c}  \; \; .
\end{equation*}
If $\tau _{\varphi }^{c}$ is completely hermitian, then $\exp it\tau
_{\varphi }^{c}$ is
a completely isometric surjection, and that is also the
case for $\tau _{\exp it\varphi }^{c}$.
    From
Corollary
\ref{unit},
$\psi (t)=\exp it\varphi $ is a unitary element of $A_{\ell}(X).$
Since $t\mapsto \psi (t)$ is a norm continuous one-parameter group of
unitaries in the $C^{*}$-algebra $A_{\ell }(X),$ it follows that $\varphi $
is a self-adjoint element in $A_{\ell }(X).$
\endproof

\section{Some applications}\label{applications}

If $x, y \in B(H)$ then we say that $x \perp y$ if
$x^* y = 0$.  Similarly for subsets $E , F \subseteq B(H)$,
we write $E \perp F$ if $x^* y = 0$ for all
$x \in E, y \in F$.

\begin{theorem} \label{tot2}  If $P$ is a projection on an operator space
$X$, then
the following are equivalent (and are also equivalent
to the conditions in Proposition \ref{Mproject}):
\begin{itemize}
\item [{\rm (a)}]  $P$ is a complete left $M$-projection.
\item [{\rm (b)}]   $\tau^c_P$  is completely contractive.
\item [{\rm (c)}]  $P$ is an orthogonal  projection in the
$C^*$-algebra $A_{\ell}(X)$.
\item [{\rm (d)}]  $P \in M_\ell(X)$ with multiplier norm $\leq 1$.
\item [{\rm (e)}]  There exists  an embedding
$\sigma : X \hookrightarrow B(H)$ such that
$$\sigma(P(X)) \perp \sigma((I-P)(X)) \; . $$
\end{itemize}
In {\rm (e)}, $\sigma$ may be taken to be the Shilov embedding.
\end{theorem}

\proof
That (b) is equivalent to (d), and that (a) implies
(d) follows from Theorem \ref{lmu}
and Proposition \ref{Mproject}.
Let us assume (d). If we use a
Shilov embedding $\sigma_0:X\hookrightarrow B(K,H)$, it follows that
$P=L^{\sigma_0}(b)$, where
$b\in B(H)$ is a contraction.
Since $L^{\sigma_0}$ is one-to-one, $b^{2}=b$,
and from elementary operator theory, $b=b^{*}$ is an orthogonal
projection on $H$.
Thus $P=L^{\sigma_0}(b)$ is an orthogonal projection in $A_{\ell}(X)$ and we
have (c). Given (c), it
is immediate that $P$ is the image of an orthogonal  projection in
$A_\ell^{\sigma_{0}}(X)$.  From Proposition
\ref{Mproject} that implies (a).

Given (a), there exists by Proposition \ref{Mproject}
an embedding
$\sigma : X \hookrightarrow B(H)$, and an orthogonal  projection
$e \in B(H)$ with $\sigma(Px) = e \sigma(x)$ for all
$x \in X$.   From the above discussion we see that we
can take $\sigma$
to be a Shilov embedding and it is evident that (e) holds
for this $\sigma$.
Finally, given (e), we will show that $P$ is adjointable in
the sense of \S 4 of
\cite{BlecherShilov}.  If $x, y \in X$, then
$$\sigma(Px)^* \sigma(y) = \sigma(Px)^* \sigma(Py + (I-P)y)
= \sigma(Px)^* \sigma(Py) \; , $$
and also
$$\sigma(x)^* \sigma(Py)  = (\sigma(Py)^* \sigma(x))^* =
(\sigma(Py)^* \sigma(Px))^* = \sigma(Px)^* \sigma(Py) \; . $$
Since these are equal, $P$ is adjointable. It follows from
\cite{BlecherShilov}
that $P$ satisfies (c).
\endproof

We now wish to
investigate the $C^*$-algebra $A_\ell(X)$ in the case that
$X$ is the operator space dual of an operator space.

The following is well-known (see
Lemma A.4.2 in \cite{EffrosRuanbook} and \cite{BonsallDuncan}, I.10.10).

\begin{lemma}
\label{selfadjlem}Given an operator $d$ on a Hilbert space $H$ with $\left\|
d\right\| \leq 1,$ we have that $d=d^{*}$ if and only if $\left\|
1+itd\right\| \leq \sqrt{1+t^{2}}$ for all $t\in \Bbb{R}$. If $d$ is
an element of
a unital Banach algebra $A$ such that $\left\| 1+itd\right\| \leq
\sqrt{1+t^{2}}$ for all
$t\in {\Bbb R}$, then it is hermitian in $A$.
\end{lemma}

\begin{lemma}\label{contractive}
Given an operator space $X$ and a left multiplier $\varphi:X \to X$ such
that $\Vert \varphi \Vert_{M_{\ell}(X)}\leq 1$,
it follows that $\tau_{\varphi}^{c}$ is
a left multiplier of
$C_{2}(X)$ with $\Vert \tau_{\varphi}^{c} \Vert_{M_{\ell}(C_{2}(X))}\leq 1$. If
$\varphi$ is a self-adjoint or adjointable left multiplier, then the
same is true for
$\tau_{\varphi}^{c}$.
\end{lemma}
\proof Let us suppose that $X\subseteq B(H_{0})$ is a Shilov
embedding and that $\varphi(x)=bx$, where $b\in B(H_{0})$. We have a natural
embedding $\sigma:C_{2}(X)\hookrightarrow B(H_{0}^{2})$ defined by
\[\sigma\left (\left[
\begin{array}{c}
x \\
y
\end{array}
\right]\right )
=\left[
\begin{array}{cc}
x & 0 \\
y & 0
\end{array}
\right]
\]
and we have that
\[\sigma\left (\tau^c_{\varphi}\left (\left[
\begin{array}{c}
x \\
y
\end{array}
\right]\right )\right )
=\left[
\begin{array}{cc}
bx & 0 \\
y & 0
\end{array}
\right]=\left[
\begin{array}{cc}
b & 0 \\
0 & I
\end{array}
\right]\sigma\left (\left[
\begin{array}{c}
x \\
y
\end{array}
\right]\right )\]
where $b\oplus I$ is a contractive left multiplier of
$\sigma(C_{2}(X))$. If $b$ is
self-adjoint, then that is also the case for $b\oplus I$, and if the real and
imaginary parts of $b$ are left multipliers, that is also the case
for $b\oplus I$,
hence the remaining assertions are evident.
\endproof

\begin{theorem}
\label{dual}If $X$ is the operator space
dual of an operator space, then
$A_{\ell }(X)$ is a von Neumann algebra.
\end{theorem}
\proof

Let us suppose that $X$ is the dual of the operator space $X_{*}$. We have that
$A_{\ell }(X)$ is a Banach subalgebra of $CB(X).$ On the other hand, we
may identify $CB(X)$ with the operator space dual $(X\hat{\otimes}
X_{*})^{*}. $ To show that
$A_{\ell }(X)$ is a dual Banach space it suffices
to prove that it is closed in the weak$^{*}$ topology in $CB(X)$, and for
that it suffices to prove that its unit ball
$D = A_{\ell }(X)_{1}$ is weak$^{*}$ closed in the unit ball
$CB(X)_{1}.$
Since $X\hat{\otimes}X_{*}$ is
the norm completion of $X\otimes X_{*},$ the latter determines the same
topology on $CB(X)_{1},$ and thus given $\varphi _{\nu },\varphi \in
CB(X)_{1}$, $\varphi _{\nu }\rightarrow \varphi $ in the weak$^{*}$ topology
if and only if $\varphi _{\nu }(x)\rightarrow \varphi (x)$ in the weak$^{*}$
topology for each $x\in X.$

Suppose that
$\varphi _{\nu }\in
D_{sa}$, $\varphi \in
CB(X)_{1}$ and that $\varphi _{\nu }(x)\rightarrow \varphi (x)$ in the weak$
^{*}$ topology for each $x\in X.$
If we use the duality
\begin{equation*}
(R_{2}[X_{*}])^{*}=C_{2}(X)
\end{equation*}
it is evident that
\begin{equation*}
\tau _{\varphi _{\nu }}^{c}\left( \left[
\begin{array}{l}
x \\
y
\end{array}
\right] \right) =\left[
\begin{array}{c}
\varphi _{\nu }(x) \\
y
\end{array}
\right] \rightarrow \left[
\begin{array}{c}
\varphi (x) \\
y
\end{array}
\right]
\end{equation*}
in the weak$^{*}$ topology for any $x,y\in X.$
Hence $\tau^c_{\varphi_\nu}  \rightarrow \tau^c_{\varphi}$
in the weak*-topology of
$CB(C_2(X)) \cong (C_2(X) \hat{\otimes} R_{2}[X_{*}])^{*}$.
This follows by considerations similar to those mentioned
at the end of the last paragraph, but with
$X$ replaced by $C_2(X)$. From Lemma \ref{selfadjlem} and Lemma
\ref{contractive}, $\Vert 1 + it
\tau^c_{\varphi_\nu}
\Vert_{cb} \leq
\sqrt{1+ t^2}$ for all $t \in {\Bbb R}$, and since norm closed balls are
weak*-closed, $\Vert 1+it\tau^c_{\varphi}
\Vert_{cb} \leq
\sqrt{1+ t^2}$ for all $t \in {\Bbb R}$. From Lemma \ref{selfadjlem},
$\tau^c_{\varphi}$
is a hermitian element of the Banach algebra $CB(C_{2}(X))$, and we
have from Theorem
\ref{cch} that $\varphi \in A_{\ell}(X)_{sa}$. On the other hand, since
the norm closed unit
balls in $CB(X)$ are weak$^{*}$ closed, $\varphi\in D_{sa}$. We conclude that
$D_{sa}$ and $A_{\ell}(X)_{sa}$ are weak$^{*}$ closed.

Finally, let us suppose that
$\varphi _{\nu }\in D,$ $\varphi
\in CB(X)_{1},$ and that $\varphi _{\nu }\rightarrow \varphi $ in
the weak$^{*}$ topology.
Since $CB(X)_{1}$ is compact in the weak$^{*}$
topology, by passing to a subnet twice we may assume that
\[
\mathrm{Re}
\,\varphi _{\nu }=(1/2)(\varphi_{\nu}+\varphi_{\nu}^{*})\rightarrow \psi_{1}\]
and
\[\mathrm{Im}\,\varphi _{\nu
}=(1/2i)(\varphi_{\nu}-\varphi_{\nu}^{*})\rightarrow \psi_{2}
\]
in the weak$^{*}$ topology (we are using the involution in
$A_{\ell}(X)$). It follows
that
$\varphi =\psi _{1}+i\psi _{2},$
and from the previous argument, $\psi _{i}\in
A_{\ell }(X)_{sa}.$ As in the self-adjoint case we have that
$\Vert\varphi\Vert \leq 1$
hence $D$ and
therefore
$A_{\ell }(X)$ are weak$^{*}$ closed.
\endproof
 
This result is an important tool in our theory, since it allows the
introduction of von Neumann algebra methods to the study of dual
operator spaces.   For example, we see immediately that a dual
operator space $X$ has no nontrivial complete left $M$-projections
if and only if $A_\ell(X) =  \Bbb{C}$.  In general,
the set of complete left $M$-projections on a dual
operator space $X$ is a complete lattice;
and there is a spectral theorem
for left adjointable operators on $X$.   We plan to discuss more such
consequences in the sequel to this paper.
 
\begin{theorem}
\label{tem}
If $X$ is a dual operator space then any
$\varphi \in A_\ell(X)$ is weak*-continuous.
\end{theorem}
\proof
It suffices to prove that $\varphi$ is weak$^{*}$ continuous on the
unit ball of $X$. Since
$A_\ell(X)$ is a von Neumann algebra,
$\varphi$ is a norm limit of a sequence $\varphi_{n}$, where each
$\varphi_{n}$ is a
linear combination  of projections. The restrictions of these
mappings to the unit
ball of $X$ converge uniformly. From Theorem \ref{tot2} the
projections in $A_\ell(X)$
are the $M$-projections on $X$, and from Proposition
\ref{weakstarcont} they are weak$^{*}$ continuous. It follows that each
$\varphi_{n}$ is weak$^{*}$ continuous, and since a uniform limit of
weak$^{*}$ continuous
functions is weak$^{*}$ continuous, we conclude that $\varphi$ is
weak$^{*}$ continuous on
the unit ball of $X$.
\endproof

We note that we can also prove the above corollary by using the fact
that any element
$a\in A_\ell(X)$ with $0\leq a \leq 1$ is a corner of a projection in
$M_{2}(A_{\ell}(X))=A_{\ell}(M_{2}(X))$.

\begin{corollary}\label{automcont} If $X$ is a dual operator space and it is an
operator $A$--$B$-bimodule for
$C^{*}$-algebras $A$ and $B$, then the mapping $x\mapsto axb$ for
$a\in A$ and $b \in B$ is automatically weak$^{*}$ continuous on $X$.
\end{corollary}
\proof From \cite{CES} Corollary 3.2, there exists a completely
isometric embedding
$\Theta:X\hookrightarrow B(H)$ and $*$-representations
$\pi_{1}$ and $\pi_{2}$ of $A$ and
$B$ respectively on $H$ such that
$\Theta(axb)=\pi_{1}(a)\Theta(x)\pi_{2}(b)$ for $a\in
A$, $b\in B$ and $x\in X$. Changing notation, let us assume that
$X,A,B\subseteq
B(H)$. Since the mapping
$x\mapsto ax$ is in $A_{\ell}(X)$, we have from Theorem
\ref{tem} that it is weak$^{*}$
continuous. On the other hand since $y\mapsto yb$ is in
$A_{r}(X)$, it is also weak$^{*}$ continuous. It follows that
\[
x \mapsto axb=(ax)b
\]
is weak$^{*}$ continuous.
\endproof

\begin{corollary}\label{univdual}
\label{ndm}  Any dual operator space $X$ is a
normal dual $A_\ell(X)$--$A_r(X)\mbox{-bimod}$\-ule in the sense of
\cite{EffrosRuanbimod}, i.e, the trilinear mapping
    $A_\ell(X) \times X \times A_r(X)\rightarrow X$
is weak*-continuous in each variable.
\end{corollary}
\proof
Suppose that $a_i \in A_\ell(X)$ is a net
converging  weak* to $a$, and that $x \in X$.  Then since
the weak*-topology on $A_\ell(X)$ is inherited from
$CB(X) = (X \hat{\otimes}  X_*)^*$, we have that
$\psi(a_i (xb)) \rightarrow \psi(a (xb))$ for $\psi \in X_*$.
The same argument applies to the third variable, and continuity in
$x$ follows from the
previous corollary.
\endproof

Some further applications of these results to operator modules are
given in \cite{Blecherdual}.
 
We may also use Theorem \ref{multiplierth} to study functorial properties of
left multiplier mappings.
Given a subspace $Y$ of an operator space $X$ and
a left multiplier mapping $\varphi :X\rightarrow X$ such that $\varphi
(Y)\subseteq Y,$ it is trivial that the restriction $\varphi ^{\prime
}=\varphi _{|Y}$ is a left multiplier of $Y$. The following is perhaps less
evident.

\begin{proposition}
Suppose that $Y$ is a closed subspace of an operator space $X$, and that $
\varphi \in
M_{\ell }(X)_1$ is such that
$\varphi (Y)\subseteq Y.$ Then the induced quotient mapping $\varphi
^{\prime \prime }:X/Y\rightarrow X/Y$ is an element of
$M_{\ell }(X/Y)_1$.   If
in addition $\varphi \in A_{\ell }(X)$, and $\varphi ^{*}(Y)\subseteq Y,$
then $\varphi ^{\prime} \in  A_{\ell }(Y)$ and
$\varphi ^{\prime \prime }\in A_{\ell }(X/Y).$
\end{proposition}
\proof
Let us suppose that $\left\| \tau _{\varphi }^{c}\right\| _{cb}\leq 1.$ From
the definition of the quotient operator space structure (applied to
rectangular matrices) we have the identification
\begin{equation*}
C_{2}(X/Y)=C_{2}(X)/C_{2}(Y)
\end{equation*}
Thus an element $\left[
\begin{array}{l}
\bar{x} \\
\bar{y}
\end{array}
\right] \in C_{2}(X/Y),$ with norm less than 1 is the quotient image of an
element $\left[
\begin{array}{l}
x \\
y
\end{array}
\right] \in C_{2}(X)$   with norm less than 1. We have that $\left[
\begin{array}{c}
\varphi ^{\prime \prime }(\bar{x}) \\
\bar{y}
\end{array}
\right] $ is the quotient image of  $\left[
\begin{array}{c}
\varphi (x) \\
y
\end{array}
\right] $, and thus
\begin{equation*}
\left\| \left[
\begin{array}{c}
\varphi ^{\prime \prime }(\bar{x}) \\
\bar{y}
\end{array}
\right] \right\| \leq \left\| \left[
\begin{array}{c}
\varphi (x) \\
y
\end{array}
\right] \right\| <1,
\end{equation*}
from which it follows that
$\left\| \tau _{\varphi^{\prime \prime }}^{c}\right\| \leq 1.$ A
similar argument can be used on matrices.
 
Let
$B$ denote the $\varphi \in A_{\ell }(X)$ such that
$\varphi (Y)\subseteq Y$ and $\varphi ^{*}(Y)\subseteq
Y$.  Then
$B$ is a *-subalgebra of $A_{\ell }(X)$,
and  we have from
above that $\varphi \mapsto \varphi ^{\prime \prime }$ is a norm decreasing
unital homomorphism from
$B$ into $M_{\ell }(X/Y).$ It is
evident from Lemma \ref{selfadjlem} that the image of a self-adjoint element
of $B$ is again self-adjoint, and thus this mapping sends
$B = B_{sa} + iB_{sa}$
into $A_{\ell }(X/Y).$%
\endproof

\section{Examples}
\label{examples}

\noindent {\bf 6.1.}  It was shown in \S 4.22 of
\cite{BlecherShilov}
that if $X$ is a Banach space, then
$A_\ell(MIN(X))$ coincides with
the classical centralizer algebra $Z(X)$ of $X$. Since the
projections of $Z(X)$ are
the $M$-projections, whereas the projections in
$A_\ell(MIN(X))$ are the complete left $M$-projections on the
operator space $MIN(X)$, these mappings coincide. It follows that the
complete right $M$-summands and complete right
$M$-ideals of $MIN(X)$ are the
$M$-summands and $M$-ideals of $X$. Since in general
$MAX(X)^{*}=MIN(X^{*})$, we also
see that complete left $L$-projections (respectively, complete
right $L$-summands) in $MAX(X)$ are the $L$-projections (respectively,
$L$-summands) in $X$.
\vspace{5 mm}

\noindent {\bf 6.2.}  From Lemma \ref{2sid}, the ``complete
$M$-projections'' considered in
\cite{ERcmp} are just the complete left $M$-projections which are
also complete
right $M$-projections.   Hence it follows from \ref{list} (a) that
the ``complete
$M$-summands'' coincide with the complete left $M$-summands which are
also complete
right $M$-summands. In turn, the ``complete $M$-ideals'' of
\cite{ERcmp}  are the complete left
$M$-ideals which are also complete right $M$-ideals.
There is an operator space version of the centralizer algebra of a
Banach space which is appropriate to this
``complete two-sided'' theory, which we will consider
elsewhere.  One description of this algebra is the left adjointable
multipliers which are also right adjointable.

\vspace{5 mm}

    As we indicated in the
introduction, the complete right $M$-ideals in a
$C^*$-algebra coincide with the
closed right ideals, and the complete right $M$-summands are the
``principal right ideals'' of the form $eA$ for an orthogonal
projection $e \in M(A)$,
the  multiplier algebra of $A$.  (Indeed in \cite{BESZ}
    we show that the word ``complete'' is unnecessary here).  We
consider two generalizations
of this observation.
 
\vspace{5 mm}
 
\noindent {\bf 6.3.}
If $A$ is an operator algebra, we let $LM(A)$ be the
left multiplier algebra of $A$ (this is
equal to $A$ if $A$ is unital).\addtocounter{claim}{3}
\begin{proposition}
If $A$ is a (possibly
non-self-adjoint)
    operator algebra with
contractive approximate identity, then the complete right $M$-summands of
$A$ are exactly the principal
    right ideals $eA$ for an orthogonal  projection $e \in LM(A)$.
The complete right $M$-ideals  of
$A$ are exactly the closed right ideals of $A$ which possess a left
contractive approximate identity.
\end{proposition}
\begin{proof}
    The  first assertion is a
consequence of 4.17 in \cite{BlecherShilov},  which states that
$M_{\ell}(A) = LM(A)$
(this fact may also be proved more directly).
Hence the complete left $M$-projections
on $A$ are exactly the  orthogonal projections $e \in LM(A)$.
If $A$ is unital, this part of the argument would be easier.

It is well known that $A^{**}$ is an operator algebra with the Arens
product.
If $J$ is a complete right $M$-ideal of $A$, then $J^{**} = J^{\perp \perp}
= \bar{J}^{w*}$ is, by the first part, equal to a principal right
ideal $eA^{**}$.  Here $e \in A^{**}$ is an orthogonal  projection.
Considered as subsets of $A^{**}$, we have $J A \subset J^{**}$.  But also
$J A \subset A$.  So $J A \subset J^{**} \cap A = J$ by basic functional
analysis.  So $J$ is a right ideal of $A$.
Since $A^{**}$ is unital, $e \in J^{**}$, and
$e$ is a left identity for $J^{**}$.
There exists a net in $Ball(J)$ which converges to
$e$ in the weak* topology.  By a well known argument
using the fact that the weak closure of a convex set equals its
norm closure, one may replace the above net with a
left contractive approximate identity for $J$
(see e.g. Theorem 2.2 in \cite{ERns} for details).

Conversely, if $J$ is a  closed right ideal of $A$ with a
contractive left approximate identity, then $ J^{**}$ is a
subalgebra of $A^{**}$ with a left identity $e$ of norm $1$
by
   \cite{BonsallDuncan} 28.7.  Note that $e$ is an
    orthogonal projection
in $A^{**}$.   Moreover $J^{**} A^{**} \subset  J^{**}$
by a routine argument
approximating elements in $X^{**}$ by weak*-converging nets
of elements in $X$ (see e.g. \cite{ERns}  Theorem 2.2).
    We have
$$J^{**} = e J^{**}  \subset e A^{**}  \subset J^{**} \; \; .$$
Thus $J^{**} = e A^{**}$ is a complete right $M$-summand of $A^{**}$ by the
first part, so that $J$ is a complete right $M$-ideal of $A$.\end{proof}

There is a stronger result  due to
Zarikian \cite{Zarikianthesis} which is valid, in which
the hypothesis ``complete'' is weakened.
 
\vspace{5 mm}

\noindent {\bf 6.5.}
We next consider the one-sided $M$-structure of Hilbert $C^*$-modules.
\addtocounter{claim}{1}
\begin{theorem}The complete right $M$-ideals in a right Hilbert
$C^*$-module are exactly the closed right submodules.
The complete right $M$-summands
are the orthogonally complemented
right submodules.
\end{theorem}
\proof
The last statement is clear from Theorem
\ref{tot2}, since for any right Hilbert
$C^*$-module $X$, we have that $A_\ell(X)$ is the algebra of
adjointable operators on $X$.  Thus
the complete left $M$-projections are exactly the
adjointable projections on $X$.
 
We may assume by Cohen's factorization theorem
that the right Hilbert
$C^*$-module $X$ is full over a $C^*$-algebra $D$.
We refer to \cite{Paschke} for information on
self-dual $W^*$-modules.
We will also use the following
facts mentioned at the end of \S 5 in \cite{BlecherShilov}.
We believe that these facts are essentially folklore.  Namely,
the second dual of the linking $C^*$-algebra for $X$,
is the ``linking $W^*$-algebra'' for $X^{**}$, and the last
space $X^{**}$
is a self-dual right $C^*$-module over $D^{**}$.  As is
often very helpful in  $C^*$-module theory, one may
    view the computations below as taking place within
these linking algebras.
 
If $Y$ is a complete
right $M$-ideal of a full right Hilbert
    $C^*$-module
$X$ over $D$, then $Y^{\perp \perp} = \bar{Y}^{w*} = Y^{**}$
is a  complete
right  $M$-summand of $X^{**}$.  But by the above,
    the complete left $M$-projection on  $X^{**}$ corresponding to
$Y^{\perp \perp}$ is a $D^{**}$-module map.  Thus
$Y^{\perp \perp} = Y^{**}$ is a $D^{**}$-submodule
of $X^{**}$.  Hence
viewed as subsets of $X^{**}$, we have that
$Y D \subset
Y^{\perp \perp} \cap X$.  But the latter
space is just $Y$,  by basic functional analysis.
Thus $Y$ is a $D$-submodule of $X$.

Conversely, if
$Y$ is a $D$-submodule of a full Hilbert $C^*$-module
$X$ over $D$,
then $Y^{\perp \perp} =  \bar{Y}^{w*} = Y^{**}$ is a weak*-closed
$D^{**}$-submodule
of $X^{**}$.  This may be seen from a
routine argument
approximating elements in $X^{**}$ by weak*-converging nets
of elements in $X$.

It is a well-known fact
that a weak*-closed submodule of a self-dual $W^*$-module
is orthogonally complemented.  Since we
are not aware of a precise reference in the literature
for this we give a short proof:
Suppose that $Z$ is a weak* closed submodule of a self-dual right
$C^*$-module $X$ over a von Neumann algebra $M$.  Suppose that
$N$ is the von Neumann algebra acting on the left of $X$ (which may
be viewed as the set of bounded adjointable $M$-module maps
on $X$), and let ${\mathcal I}$ be the weak*-closure of
$Z \bar{X}$
in $N$.   Here $Z \bar{X}$ is the span
of the rank one operators
$z \otimes x$ for $z \in Z, x \in X$.
   Let ${\mathcal L}$ be the linking $W^*$-algebra of $X$, and
consider the
subspace of  ${\mathcal L}$ which has ${\mathcal I}$ and $Z$ as
its first row, and zero entries on the second.  This subspace is
a weak* closed right ideal in ${\mathcal L}$, and therefore
equals $E {\mathcal L}$ for a projection $E \in {\mathcal L}$.
It is easy to check that $E$ has only one nonzero entry, namely
its 1-1-entry, and this is the desired projection in $N$ onto
$Z$.

Since $X^{**}$ is a self-dual $W^*$-module it follows that
$Y^{\perp \perp}$ is complemented,
i.e. there exists an adjointable projection
on $X^{**}$ with range
$Y^{\perp \perp}$.  So by the first part,
$ Y^{\perp \perp}$ is
a complete right  $M$-summand, so that $Y$ is a
complete right  $M$-ideal.
\endproof
 
The following apparently new result follows from this and Theorem
\ref{list} (it can
also be proved directly).
 
\begin{corollary}  There is at most one contractive linear
projection from a Hilbert $C^*$-module onto a closed submodule.
If there exists such a projection then it is an
adjointable module map,
so that the $C^*$-submodule is complemented orthogonally in
the sense of  Hilbert $C^*$-module theory.
\end{corollary}
    \vspace{5 mm}

\noindent{\bf 6.8.}
We may also describe the one-sided $M$-structure
of various Hilbertian operator spaces. In this discussion we let $H$
denote a Hilbert
space, and $H_{c}$, $H_{r}$ and $H_{0}$ denote the column, row, and
Pisier's self-dual
quantizations of $H$ (see, e.g. \cite{EffrosRuanbook}).
\addtocounter{claim}{1}
\begin{lemma}\label{Lvarious}
Let $\xi, \eta \in H$. Then
\begin{equation} \label{1st}
           \left\|
           \begin{bmatrix}
                   \xi\\
                   \eta\\
           \end{bmatrix}
           \right\|_{C_2(H_c)}  \!\!\!= \left\|
        \begin{bmatrix}
                   \xi\\
                   \eta\\
           \end{bmatrix}
           \right\|_{C_2[H_c]}  \!\!\!= \left\|
           \begin{bmatrix}
                   \xi \!\!\!& \eta
           \end{bmatrix}
           \right\|_{R_2(H_r)} = \left\|
           \begin{bmatrix}
                   \xi \!\!\! & \eta
           \end{bmatrix}
           \right\|_{R_2[H_r]} = \!\!\! \; \; \sqrt{\|\xi\|^2 + \|\eta\|^2}.
\end{equation}
This also coincides with the $C_2(MAX(H))$
and the $C_2[MIN(H)]$
norms.
If, in addition, $\xi \perp \eta$, then
\begin{equation} \label{2nd}
           \left\|
           \begin{bmatrix}
                   \xi\\
                   \eta
           \end{bmatrix}
           \right\|_{C_2(H_r)}
=  \left\|
           \begin{bmatrix}
                   \xi\\
                   \eta
           \end{bmatrix}
           \right\|_{C_2(MIN(H))}
= \left\|
           \begin{bmatrix}
                   \xi  \!\!\! & \eta
           \end{bmatrix}
           \right\|_{R_2(H_c)}
= \max\{\|\xi\|, \|\eta\|\},
\end{equation}
\begin{equation} \label{3rd}
           \left\|
           \begin{bmatrix}
                   \xi\\
                   \eta
           \end{bmatrix}
           \right\|_{C_2[H_r]} =
     \left\|
           \begin{bmatrix}
                   \xi\\
                   \eta
           \end{bmatrix}
           \right\|_{C_2[MAX(H)]} =
\left\|
           \begin{bmatrix}
                   \xi \!\!\! & \eta
           \end{bmatrix}
           \right\|_{R_2[H_c]} = \|\xi\| + \|\eta\|,
\end{equation}

\begin{equation} \label{4th}
           \left\|
           \begin{bmatrix}
                   \xi\\
                   \eta
           \end{bmatrix}
           \right\|_{C_2(H_o)} =
\sqrt[4]{\|\xi\|^4 + \|\eta\|^4} \; \; \; , \;
\; \text{and} \; \; \;
           \left\|
           \begin{bmatrix}
                   \xi\\
                   \eta
           \end{bmatrix}
           \right\|_{C_2[H_o]} =
(\|\xi\|^{\frac{4}{3}} + \|\eta\|^{\frac{4}{3}})^{\frac{3}{4}} \; \; .
\end{equation}
\end{lemma}
\proof We may assume that neither
$\xi$ nor $\eta$ is zero. Equations (\ref{1st}) and
(\ref{2nd}) are well known,
for example (\ref{1st})
follows from the completely isometric identifications $C_2(H_c) \cong
(H^2)_c \cong C_2[H_c]$ and $R_2(H_r) \cong (H^2)_r \cong R_2[H_r]$,
where $H^2 =
H \oplus H$.  The first assertion after (\ref{1st}) follows from (\ref{1st}),
(\ref{ineq1}),  and the fact that $MAX$ dominates
the other operator space structures. Similarly the second assertion
follows from
(\ref{1st}) and (\ref{ineq2}).
For equation (\ref{3rd}), we
will use the completely isometric
identifications $C_2[H_r] \cong
\mathcal{T}(\overline{H},\mathbb{C}^2)$ and
$R_2[H_c] \cong \mathcal{T}(\overline{\mathbb{C}^2},H)$,
where for Hilbert spaces $K$ and $L$, $\mathcal{T}(\overline{K},L)$
denotes the family of trace-class operators
from the conjugate Hilbert space of $K$ to $L$.
Under the first identification, $\begin{bmatrix} \xi\\ \eta
\end{bmatrix}$ corresponds to the mapping
$S:\overline{H} \to \mathbb{C}^2$ defined by
\[
           S(\overline{\zeta}) =
           \begin{bmatrix}
                   \langle\overline{\zeta},\overline{\xi}\rangle\\
                   \langle\overline{\zeta},\overline{\eta}\rangle
           \end{bmatrix} =
           \begin{bmatrix}
                   \langle\xi,\zeta\rangle\\
                   \langle\eta,\zeta\rangle
           \end{bmatrix}
           \text{ for all } \zeta \in H.
\]
Under the second identification,
$\begin{bmatrix} \xi & \eta \end{bmatrix}$ corresponds to the mapping
$T:\overline{\mathbb{C}^2} \to H$ defined by
\[
           S\left(
           \overline{
                   \begin{bmatrix}
                           a\\
                           b
                   \end{bmatrix}}
           \right) = \overline{a} \xi + \overline{b} \eta
           \text{ for all } a, b \in \mathbb{C}.
\]
Routine calculations then show that
\[
           \left\|
           \begin{bmatrix}
                   \xi\\
                   \eta
           \end{bmatrix}
           \right\|_{C_2[H_r]} =
    \|S\|_{\mathcal{T}(\overline{H},\mathbb{C}^2)} =
           \|\xi\| + \|\eta\|
\]
and
\[
           \left\|
           \begin{bmatrix}
                   \xi & \eta
           \end{bmatrix}
           \right\|_{R_2[H_c]} =
           \|T\|_{\mathcal{T}(\overline{\mathbb{C}^2},H)} =
           \|\xi\| + \|\eta\|.
\]
Since the norm on $C_2[MAX(H)]$ dominates
the $C_2[H_r]$ norm, it must be equal to $\|\xi\| + \|\eta\|$
too.
Finally, to prove equation (\ref{4th}), we compute
\[
           \left\|
           \begin{bmatrix}
                   \xi\\
                   \eta
           \end{bmatrix}
           \right\|_{C_2(H_o)} =
           \left\|
           \begin{bmatrix}
                   \langle\xi,\xi\rangle\\
                   \langle\xi,\eta\rangle\\
                   \langle\eta,\xi\rangle\\
                   \langle\eta,\eta\rangle
           \end{bmatrix}
           \right\|^{1/2} =
           \left\|
           \begin{bmatrix}
                   \|\xi\|^2\\
                   0\\
                   0\\
                   \|\eta\|^2
           \end{bmatrix}
           \right\|^{1/2} = \sqrt[4]{\|\xi\|^4 + \|\eta\|^4}.
\]
The second statement in (\ref{4th}) may
be
seen from the following
argument of Pisier:  it is shown in \cite{Pisier}, Theorem
2.3 that
$C_2[H_o] = H_o \otimes_h C_2$ is the midway interpolant between
$H_c \otimes_h C_2$ and $H_r \otimes_h C_2$.  But the first space
may be thought of as the Schatten 2-class, and the second space as
the trace class.  Hence $H_o \otimes_h C_2$ may be identified
with the Schatten $\frac{4}{3}$-class.    From this the claimed
statement follows easily.
\endproof

\begin{proposition}\label{table1}
The one-sided
$M$- and $L$-projections for the various
    quantizations of $H$ are given by the following table.
\begin{center}
           \begin{tabular}{| l | c | c | c | c |}
           \hline
           {\rm Quantization} & {\rm Left $M$-Proj's} & {\rm Right
$L$-Proj's} & {\rm Right
                   $M$-Proj's} & {\rm Left $L$-Proj's} \\
           \hline
           $H_c$ & $Proj(B(H))$ & $\{0, I\}$ & $\{0, I\}$ &
               $Proj(B(H))$ \\
           \hline
           $H_r$ & $\{0, I\}$ & $Proj(B(H))$ &
                   $Proj(B(H))$ & $\{0, I\}$ \\
           \hline
           $H_o$ & $\{0, I\}$ & $\{0, I\}$ & $\{0, I\}$ & $\{0, I\}$\\
           \hline
           $MIN(H)$ & $\{0, I\}$ & $Proj(B(H))$ & $\{0, I\}$ & $Proj(B(H))$\\
        \hline
           $MAX(H)$ & $Proj(B(H))$ & $\{0, I\}$ & $Proj(B(H))$ & $\{0, I\}$\\
           \hline
           \end{tabular}
\end{center}
where $Proj(B(H))$ is the family of orthogonal projections on $H$
and $I$ is the
    identity on $H$.  The first three rows of the table
are  also valid for {\em complete}
one-sided  $M$- and $L$-projections.
For $MAX(H)$ however (respectively, $MIN(H)$), there are no nontrivial
complete left- or right- $M$-projections (respectively, $L$-projections).
\end{proposition}

\proof
We begin by noting that a one-sided $M$- or $L$-projection $P$ for
any of the five
quantizations of $H$ is necessarily an orthogonal
projection on $H$ since any
    such $P$ is a contractive linear idempotent.

Our procedure for showing that a given entry in the
table is $\{0, I\}$ is to suppose the contrary.  Then there
exists an orthonormal set $\{ \xi , \eta \}$  with $P \xi = \xi$
and $P\eta = 0$.  The fact that
$\Vert \eta + \xi \Vert = \sqrt{2}$ leads to a contradiction if one appeals
to the appropriate formula in the previous lemma. For example,
our assertion that the entries in the second and third columns of the
first row
are as small as possible follows immediately from this argument and
(\ref{2nd}) and (\ref{3rd}).

The fact that the entries in the first and fourth columns of the
first row are as large as possible follows almost immediately from the previous
lemma. For example, we have from (\ref{1st}) that for any $P \in Proj(B(H))$,
\[
           \left\|
           \begin{bmatrix}
                   P\xi\\
                   (Id - P)\xi
          \end{bmatrix}
           \right\|_{C_2(H_c)} =
           \sqrt{\|P\xi\|^2 + \|(Id - P)\xi\|^2} = \|\xi\|
\]
for all $\xi \in H$ and
\[
           \|P\xi + (Id - P)\eta\| = \sqrt{\|P\xi\|^2 + \|(Id - P)\eta\|^2}
           \leq \sqrt{\|\xi\|^2 + \|\eta\|^2} = \left\|
           \begin{bmatrix}
                   \xi\\
                   \eta
           \end{bmatrix}
           \right\|_{C_2(H_c)}
\]
for all $\xi, \eta \in H$. In other words, $\nu_P^c:H_c \to C_2(H_c)$ is
an isometry and $\mu_P^c:C_2(H_c) \to H_c$ is a contraction.
Because of the completely isometric
identification $C_2(H_c) \cong (H^2)_c$, and because for maps
between Hilbert column spaces the norm coincides with the
completely bounded norm, we conclude
that $\nu_P^c$ and $\mu_P^c$ are in fact completely
contractive.  Consequently,  $P$ is a complete left
$M$-projection.  Since $C_2(H_c) \cong C_2[H_c]$ completely
isometrically, we see the other assertion.

Thus we have completed the first row.
The entries for the $H_r$ row follows by symmetry.

The entries in the first and third columns of the fourth row
    will be equal by symmetry, since $C_2(MIN(H)) \cong
R_2(MIN(H))$  isometrically.   Again arguing by contradiction
and (\ref{2nd}), we conclude that these entries are the trivial ones.
Similarly by symmetry the first and third columns of the third
row will be equal, and we use (\ref{4th}) to evaluate these.
Similarly the second and fourth entries of the fourth row are equal,
and we use the second statement after (\ref{1st}) to deduce
that the listed entries are correct here.
The remaining entries in the table are verified in just the same way.

Finally, we shall show that  $MIN(H)$ has no nontrivial complete right
$L$-projec\-tions  (from which the other final statements follow by duality and
symmetry). To that end, assume that $P \in Proj(B(H))$ is a
nontrivial complete right
$L$-projection for $MIN(H)$. Then there exist
orthonormal vectors $\xi, \eta \in H$ such that $P\xi = \xi$ and
$P\eta = 0$. But then
using successively (\ref{2nd}), the definition of a complete right
$L$-projection,
Lemma \ref{vlem}, and the ``row-version'' of  (\ref{ineq2}),  we have
\begin{eqnarray*}
          1
          &=& \max\{\|\xi\|,\|\eta\|\}\\
          &=& \left\|\begin{bmatrix} \xi & \eta
\end{bmatrix}\right\|_{R_2(MIN(H))}\\
          &=& \left\|\begin{bmatrix} P\xi & (Id - P)\xi & P\eta & (Id - P)\eta
\end{bmatrix}\right\|_{R_2(R_2[MIN(H)])}\\
&=& \left\|\begin{bmatrix} P\xi &
P\eta & (Id - P)\xi & (Id - P)\eta
\end{bmatrix}\right\|_{R_2[R_2(MIN(H))]}\\
&=& \left\|\begin{bmatrix} \xi & 0 & 0 & \eta
\end{bmatrix}\right\|_{R_2[R_2(MIN(H))]}\\
          &\geq& \sqrt{\left\|\begin{bmatrix} \xi & 0
\end{bmatrix}\right\|_{R_2(MIN(H))}^2
   + \left\|\begin{bmatrix} 0 & \eta \end{bmatrix}\right\|_{R_2(MIN(H))}^2}\\
          &=& \sqrt{\|\xi\|^2 + \|\eta\|^2}\\
          &=& \sqrt{2},
\end{eqnarray*}
a contradiction.
\endproof

Using this proposition, we can identify the one-sided summands and ideals
    in Hilbert operator spaces. For example, the
(complete)
right $M$-ideals in $H_c$ are precisely the closed subspaces of $H$,
whereas the only
right $L$-summands of $H_o$ are $\{0\}$ and $H$.

\vspace{5 mm}

\noindent {\bf 6.11.}  As a final example, we note that it is proved in
\cite{BESZ} that there exist no nontrivial
complete right or left $L$-projections on a $C^*$-algebra.  Equivalently,
    there exist no nontrivial
complete right or left $M$-projections on the
    predual of a von Neumann algebra.
In these results we may replace the word ``complete'' by
``strong'' (see \S 3).
 
\vspace{5 mm}

{\bf Remarks added March 2001:}  V. Paulsen has found an elegant
proof of Theorem 4.6 based on a 3$\times $3 matrix argument.

As we have indicated elsewhere \cite{Blecherdual}, but which
is appropriate to state here,
Theorem 4.6 facilitates
a deeper understanding of the interplay between the multiplication
operation and
the metric structure of an operator algebra.
On the one hand, it gives more
or less immediately the `BRS' characterization
of operator algebras, and the `CES' characterization of
operator modules (or more generally, the `oplication theorem'
of \cite{BlecherShilov}).   This was independently observed by Paulsen.
On the other hand it enables one to recover the
multiplication operation on a unital operator alebra from its
underlying operator
space structure.

To illustrate the second assertion, let us suppose that $A$ is an operator
algebra with an identity of norm
$1$, but that we have `forgotten' the multiplication operation on
$A$.  Let us assume
for a moment that we do `remember' the identity element $e$.
Form $M_\ell(A)$ using Theorem 4.6,
and define $\theta : M_\ell(A) \rightarrow A$ by
$\theta(T) = T(e)$.  Then the product on $A$ is given by
$ab = \theta(\theta^{-1}(a) \theta^{-1}(b))$.

If we have also forgotten the specific identity element $e$, then
we can only retrieve the product on $A$ up to a unitary $u$ with
$u, u^{-1} \in A$.   Such unitaries form a group.
Indeed they are characterized
by the Banach-Stone theorem for operator algebras (see e.g. the
last page of \cite{BlecherShilov}, or
\cite{Kadison} for the $C^*$-algebra case) as
the elements $x_0$ with the
property that the map $\pi : T \mapsto T(x_0)$ is
a completely isometric surjection
$M_{\ell}(A) \rightarrow A$.
If $A$ is a $C^*$-algebra one only needs this to
be an isometry, by Kadison's result \cite{Kadison}.
We remark that from Lemma 4.5 the unitaries in a $C^*\mbox{-algebra}$
correspond
to linear $\varphi : A \rightarrow A$ such that
$\tau^c_\varphi$ is a surjective isometry.
However in this case there are other Banach space
characterizations of unitaries - C. Akemann and N. Weaver
have shown us one
such \cite{AW}.   Given such an $x_0$ and $\pi$, we may again
recover the product as
$ab = \pi(\pi^{-1}(a)\pi^{-1}(b))$.  This is the operator
algebra product on $A$ which has this unitary
as the identity.  This is all fairly easy to see from the
Banach-Stone theorem mentioned above
and basic facts about the left multiplier
algebra.


\begin{thebibliography}{99}
\bibitem{AW} C. Akemann and N. Weaver,
Geometric characterization of some classes of operaors in
$C^*$-algebras, Preprint (2001).

\bibitem{Alfsen}  E. M. Alfsen, {\em $M$-structure and intersection
properties of balls in Banach spaces,} Israel J. Math.
   {\bf 13} (1972), 235-245.
 
\bibitem{Alfsenbook}   E. M. Alfsen, {\em  Compact convex sets and
boundary integrals,}
Erg. Math. 57,
Springer Verlag, Berlin, 1971.
 
\bibitem{AlfsenAndersen}  E. M. Alfsen and  T. Andersen,
{\em  Split faces of compact convex sets,}
    Proc.  London Math. Soc. (3) {\bf 21} (1970), 415-442.
 
 
\bibitem{AlfsenEffros}  E. M. Alfsen and E. G. Effros, {\em Structure in
real Banach spaces I \& II},  Ann. of Math. {\bf 96} (1972), 98-173.
 
 
\bibitem{AlfsenSchultz2}  E. M. Alfsen and F. Schultz, {\em On
orientation and dynamics
in operator
algebras I}, Comm. Math. Phys.(1) {\bf 194} (1998), 87-108.
 

\bibitem{Arv1}  W. B. Arveson,
{\em Subalgebras of }$C^{*}-${\em algebras,}
    Acta Math. {\bf 123 }(1969), 141-224.


\bibitem{Arv2}   W. B. Arveson,
{\em Subalgebras of }$C^{*}-${\em algebras II,}
   Acta Math. {\bf 128} (1972), 271-308.

 
\bibitem{Behrends}  E. Behrends, {\em $M$-structure and the Banach-Stone
theorem,}  Lecture Notes in Math. 736, Springer-Verlag, Berlin
(1979).
 
\bibitem{BlecherShilov}  D. P. Blecher, {\em
    The Shilov boundary of an operator space and
the characterization theorems,} J. Fnct. Anal.,
to appear (Revision of October 2000).
 
\bibitem{Blecherdual}  D. P. Blecher, {\em Multipliers and dual
operator algebras,}  J. Fnct. Anal., to appear.
 
\bibitem{BESZ}   D. P. Blecher, R. R. Smith, and
V. Zarikian, {\em One-sided projections on C$^*-$algebras},
Preprint (2002).
 
\bibitem{BlecherPaulsen}    D. P. Blecher and V. I. Paulsen,
{\em Multipliers of operator spaces
and the injective
envelope,}  Pacific Journal of Math, to appear.

\bibitem{BonsallDuncan}  F. F. Bonsall and J. Duncan,
{\em Complete normed
algebras,}      Springer-Verlag, New York-Heidelberg
(1973).

 
\bibitem{CES} E. Christensen,
E. Effros, and A. Sinclair, {\em Completely b\!\!ounded
multilinear maps and $C^{*}\mbox{-algebraic}$ cohomology}, Inventiones Mat.
{\bf 90} (1987), 279-290.
 
\bibitem{Effros} E. G. Effros, {\em Order ideals in a
C$^*-$algebra and its dual,} Duke Math. J. {\bf 30}
(1963), 391-412.
 
 
\bibitem{EffrosRuanbimod}   E. G.
Effros and Z. J. Ruan, {\em Representations of
operator bimodules and their
applications, }J. Operator Theory {\bf 19 }%
(1988), 137-157.

\bibitem{ERns}  E. G. Effros and Z. J. Ruan,
{\em On non-self-adjoint operator
algebras},  Proc. Amer. Math. Soc. {\bf 110} (1990),
915-922.

 
\bibitem{ERcmp}   E. G. Effros and Z. J. Ruan,
{\em Mapping spaces and liftings
for operator spaces},  Proc. London Math. Soc.
    {\bf 69} (1994), 171-197.
 

\bibitem{EffrosRuanbook}   E. G.
Effros and Z. J. Ruan, {\em  Operator Spaces}, Oxford University
Press, Oxford (2000).

\bibitem{EOR}   E. G. Effros, N. Ozawa, and Z. J. Ruan,
    {\em On injectivity and
nuclearity for operator spaces,} Duke Math. J., to appear.


\bibitem{Ham1}  M. Hamana, {\em    Injective
envelopes of operator systems,}
Publ. R.I.M.S. Kyoto Univ. {\bf 15} (1979), 773-785.

     \bibitem{Hamana}  M. Hamana, {\em Triple  envelopes and Silov
boundaries of operator spaces,}  Math. J. Toyama University
{\bf 22} (1999), 77-93.
 
\bibitem{HWW} P. Harmand, D. Werner and W. Werner,
{\em $M$-ideals in Banach spaces and Banach algebras,}
    Lecture Notes in Mathematics  1547,
    Springer-Verlag, Berlin - New York (1993).
 
 
\bibitem{Kadison}  R. V. Kadison,
{\em  Isometries of operator algebras,}
    Ann. of Math. (2) {\bf 54} (1951), 325-338.

\bibitem{Kirchberg} E. Kirchberg, {\em  On restricted perturbations in
inverse images and a
description of normalizer algebras in $C^{*}$-algebras}, J. Fnal Anal.
{\bf 129} (1995),
35-63.
 
\bibitem{Lance}     E. C. Lance,
{\em Hilbert C$^*$-modules - A toolkit for operator
algebraists,} London Math. Soc. Lecture Notes,
Cambridge University Press (1995).
 
\bibitem{Paschke} W. Paschke, {\em Inner product modules over
B$^*$-algebras,}
Trans. Amer. Math. Soc. {\bf 182} (1973), 443-468.

\bibitem{Pisier}  G. Pisier, {\em The operator Hilbert space
OH, complex interpolation and tensor norms,}
Memoirs Amer. Math. Soc. {\bf 585} (1996).
 
\bibitem{Prosser} R. T. Prosser, {\em On the
ideal structure of operator algebras},
Memoirs Amer. Math. Soc {\bf 45} (1963).

\bibitem{Rieffel} M. A. Rieffel,
{\em Morita equivalence for }$C^{*}-${\em algebras and
W}$^{*}-${\em algebras,} J. Pure Appl. Algebra {\bf 5 }(1974),
51-96.

\bibitem{Rieffelsu}
M. A. Rieffel, {\em Morita equivalence for operator
algebras, }Proceedings of Symposia in Pure Mathematics
    {\bf 38} Part 1
(1982), 285-298.

 
\bibitem{Ruanoper}  Z. J. Ruan, {\em Subspaces
of }$C^{*}-${\em algebras,}
J.  Fnal Anal. {\bf 76} (1988), 217-230.
 
\bibitem{Ruaninjec}  Z. J. Ruan, {\em Injectivity of operator spaces,}
Trans. Amer. Math. Soc. {\bf 315} (1989), 89-104.
 
 
\bibitem{SmithWard} R. R. Smith and J. D. Ward, {\em
$M$-ideal structure in Banach algebras,} J. Fnal
Anal. {\bf 27} (1978), 337-349.

\bibitem{Tom} M. Tomita, {\em Spectral theory of operator algebras,
II} Math. J.
Okayama Univ. {\bf 10} (1960), 19-60.
 
\bibitem{KHWerner}   K. H. Werner, {\em A
characterization of C$^*-$algebras
by nh-projections on matrix ordered spaces,}
    Preprint, Universitat des Saarlandes, 1978.

 
\bibitem{WWerner}  W. Werner, {\em Small K-groups for operator systems,}
Preprint, 1999.
 
\bibitem{Wittstock}   G. Wittstock, {\em Extensions of
completely bounded module
morphisms}, Proceedings of conference on operator algebras and
group representations, Neptum, 238-250,
Pitman (1983).
 
\bibitem{Witt}   G. Wittstock,  {\em Matrix order and
$W^{*}$-algebras in the operational
approach to statistical physical systems,}
    Comm. Math. Phys. {\bf 74} (1980), 61-70.
 
\bibitem{Zarikianthesis}  V. Zarikian, Thesis UCLA.
 
 
    \end{thebibliography}
\end{document}